\input amstex
\documentstyle{amsppt}
%----------------------------------------------------------------
% Title:     A note on Dirac spinors in a non-flat space-time of 
%            general relativity.
% Authors:   Ruslan Sharipov
% Comments:  AmSTeX, 22 pages, amsppt style
% MSC-class: 53C27, 83C60
%----------------------------------------------------------------
%           Replacement for output macro definition
%
\catcode`@=11
\redefine\output@{%
  \def\break{\penalty-\@M}\let\par\endgraf
  \ifodd\pageno\global\hoffset=105pt\else\global\hoffset=8pt\fi  
  \shipout\vbox{%
    \ifplain@
      \let\makeheadline\relax \let\makefootline\relax
    \else
      \iffirstpage@ \global\firstpage@false
        \let\rightheadline\frheadline
        \let\leftheadline\flheadline
      \else
        \ifrunheads@ %\let\makefootline\relax
        \else \let\makeheadline\relax
        \fi
      \fi
    \fi
    \makeheadline \pagebody \makefootline}%
  \advancepageno \ifnum\outputpenalty>-\@MM\else\dosupereject\fi
}
\catcode`\@=\active
%----------------------------------------------------------------
\nopagenumbers
\def\negskp{\hskip -2pt}
\def\MatGrGL{\operatorname{GL}}
\def\MatGrSL{\operatorname{SL}}
\def\MatGrSO{\operatorname{SO}}
\def\MatGrO{\operatorname{O}}
\def\GrPin{\operatorname{Pin}}
\def\GrSpin{\operatorname{Spin}}
\def\idop{\operatorname{\bold{id}}}
\def\compos{\,\raise 1pt\hbox{$\sssize\circ$} \,}
\def\vtrule{\vrule height 12pt depth 6pt}
\def\vttrule{\vrule height 12pt depth 6pt\vrule height 12pt depth 6pt}
\def\vtttrule{\vrule height 12pt depth 19pt}
\def\boxit#1#2{\vcenter{\hsize=122pt\offinterlineskip\hrule
  \line{\vtttrule\hss\vtop{\hsize=120pt\centerline{#1}\vskip 5pt
  \centerline{#2}}\hss\vtttrule}\hrule}}
\def\msum#1{\operatornamewithlimits{\sum^#1\!{\ssize\ldots}\!\sum^#1}}
\accentedsymbol\bd{\kern 2pt\bar{\kern -2pt d}}
\accentedsymbol\bPsi{\kern 1pt\overline{\kern -1pt\boldsymbol\Psi
\kern -1pt}\kern 1pt}
\def\blue#1{#1}
\catcode`#=11\def\diez{#}\catcode`#=6
\catcode`_=11\def\podcherkivanie{_}\catcode`_=8
%\catcode`~=11\def\volna{~}\catcode`~=\active
\def\mycite#1{\cite{\blue{#1}}\immediate\special{ps:
     ShrHPSdict begin /ShrBORDERthickness 0 def}}

\def\mytag#1{%
    \tag#1}
\def\mythetag#1{\thetag{\blue{#1}}\immediate\special{ps:
     ShrHPSdict begin /ShrBORDERthickness 0 def}}
\def\myrefno#1{\no#1}
\def\myhref#1#2{\blue{#2}\immediate\special{ps:
     ShrHPSdict begin /ShrBORDERthickness 0 def}}
\def\myEarXivlink{\myhref{http://arXiv.org}{http:/\negskp/arXiv.org}}
\def\myGeoCities{\myhref{http://www.geocities.com}{GeoCities}}
\def\mytheorem#1{\csname proclaim\endcsname{Theorem #1}}
\def\mythetheorem#1{\blue{#1}\immediate\special{ps:
     ShrHPSdict begin /ShrBORDERthickness 0 def}}
\def\mylemma#1{\csname proclaim\endcsname{Lemma #1}}
\def\mythelemma#1{\blue{#1}\immediate\special{ps:
     ShrHPSdict begin /ShrBORDERthickness 0 def}}
\def\mydefinition#1{\definition{Definition #1}}
\def\mythedefinition#1{\blue{#1}\immediate\special{ps:
     ShrHPSdict begin /ShrBORDERthickness 0 def}}

%----------------------------------------------------------------
\pagewidth{360pt}
\pageheight{606pt}
\topmatter
\title
A note on Dirac spinors in a non-flat space-time of 
general relativity.
\endtitle
\author
R.~A.~Sharipov
\endauthor
\address 5 Rabochaya street, 450003 Ufa, Russia\newline
\vphantom{a}\kern 12pt Cell Phone: +7-(917)-476-93-48
\endaddress
\email \vtop to 30pt{\hsize=280pt\noindent
\myhref{mailto:R\podcherkivanie Sharipov\@ic.bashedu.ru}
{R\_\hskip 1pt Sharipov\@ic.bashedu.ru}\newline
\myhref{mailto:r-sharipov\@mail.ru}
{r-sharipov\@mail.ru}\newline
\myhref{mailto:ra\podcherkivanie sharipov\@lycos.com}{ra\_\hskip 1pt
sharipov\@lycos.com}\vss}
\endemail
\urladdr
\vtop to 20pt{\hsize=280pt\noindent
\myhref{http://www.geocities.com/r-sharipov}
{http:/\negskp/www.geocities.com/r-sharipov}\newline
\myhref{http://www.freetextbooks.boom.ru/index.html}
{http:/\negskp/www.freetextbooks.boom.ru/index.html}\vss}
\endurladdr
\abstract
    Some aspects of Dirac spinors are resumed and studied in order 
to interpret mathematically the $P$ and $T$ operations in a 
gravitational field.
\endabstract
\endtopmatter
\loadbold
\loadeufb
\TagsOnRight
\document
\accentedsymbol\tbvartheta{\tilde{\overline{\boldsymbol\vartheta}}
\vphantom{\boldsymbol\vartheta}}

\rightheadtext{A note on Dirac spinors \dots}
\head
1. Algebra and geometry of two-component spinors.
\endhead
    From a mathematical point of view two-component spinors naturally arise
when one tries to understand geometrically the well-known group homomorphism
$$
\hskip -2em
\varphi\!:\MatGrSL(2,\Bbb C)\to\MatGrSO^+(1,3,\Bbb R)
\mytag{1.1}
$$
given by the following explicit formulas:
$$
\allowdisplaybreaks
\gather
\hskip -6em
\gathered
S^0_0=\frac{\overline{\goth S^1_1}\,\goth S^1_1
+\overline{\goth S^1_2}\,\goth S^1_2+\overline{\goth S^2_1}
\,\goth S^2_1+\overline{\goth S^2_2}\,\goth S^2_2}{2},\\
S^0_1=\frac{\overline{\goth S^1_1}\,\goth S^1_2
+\overline{\goth S^1_2}\,\goth S^1_1+\overline{\goth S^2_1}
\,\goth S^2_2+\overline{\goth S^2_2}\,\goth S^2_1}{2},\\
S^0_2=\frac{\overline{\goth S^1_2}\,\goth S^1_1
-\overline{\goth S^1_1}\,\goth S^1_2+\overline{\goth S^2_2}
\,\goth S^2_1-\overline{\goth S^2_1}\,\goth S^2_2}{2\,i},\\
S^0_3=\frac{\overline{\goth S^1_1}\,\goth S^1_1
-\overline{\goth S^1_2}\,\goth S^1_2+\overline{\goth S^2_1}
\,\goth S^2_1-\overline{\goth S^2_2}\,\goth S^2_2}{2},
\endgathered
\mytag{1.2}\\
\vspace{1ex}
\hskip 2em
\gathered
S^1_0=\frac{\overline{\goth S^2_1}\,\goth S^1_1
+\overline{\goth S^1_1}\,\goth S^2_1+\overline{\goth S^2_2}
\,\goth S^1_2+\overline{\goth S^1_2}\,\goth S^2_2}{2},\\
S^1_1=\frac{\overline{\goth S^2_1}\,\goth S^1_2
+\overline{\goth S^1_2}\,\goth S^2_1+\overline{\goth S^2_2}
\,\goth S^1_1+\overline{\goth S^1_1}\,\goth S^2_2}{2},\\
S^1_2=\frac{\overline{\goth S^1_2}\,\goth S^2_1
-\overline{\goth S^2_1}\,\goth S^1_2+\overline{\goth S^2_2}
\,\goth S^1_1-\overline{\goth S^1_1}\,\goth S^2_2}{2\,i},\\
S^1_3=\frac{\overline{\goth S^2_1}\,\goth S^1_1
+\overline{\goth S^1_1}\,\goth S^2_1-\overline{\goth S^2_2}
\,\goth S^1_2-\overline{\goth S^1_2}\,\goth S^2_2}{2},
\endgathered
\mytag{1.3}\\
\vspace{1ex}
\hskip -6em
\gathered
S^2_0=\frac{\overline{\goth S^1_1}\,\goth S^2_1
-\overline{\goth S^2_1}\,\goth S^1_1+\overline{\goth S^1_2}
\,\goth S^2_2-\overline{\goth S^2_2}\,\goth S^1_2}{2\,i},\\
S^2_1=\frac{\overline{\goth S^1_2}\,\goth S^2_1
-\overline{\goth S^2_1}\,\goth S^1_2+\overline{\goth S^1_1}
\,\goth S^2_2-\overline{\goth S^2_2}\,\goth S^1_1}{2\,i},\\
S^2_2=\frac{\overline{\goth S^2_2}\,\goth S^1_1
+\overline{\goth S^1_1}\,\goth S^2_2-\overline{\goth S^2_1}
\,\goth S^1_2-\overline{\goth S^1_2}\,\goth S^2_1}{2},\\
S^2_3=\frac{\overline{\goth S^1_1}\,\goth S^2_1
-\overline{\goth S^2_1}\,\goth S^1_1+\overline{\goth S^2_2}
\,\goth S^1_2-\overline{\goth S^1_2}\,\goth S^2_2}{2\,i},
\endgathered
\mytag{1.4}\\
\vspace{1ex}
\hskip 2em
\gathered
S^3_0=\frac{\overline{\goth S^1_1}\,\goth S^1_1
+\overline{\goth S^1_2}\,\goth S^1_2-\overline{\goth S^2_1}
\,\goth S^2_1-\overline{\goth S^2_2}\,\goth S^2_2}{2},\\
S^3_1=\frac{\overline{\goth S^1_1}\,\goth S^1_2
+\overline{\goth S^1_2}\,\goth S^1_1-\overline{\goth S^2_1}
\,\goth S^2_2-\overline{\goth S^2_2}\,\goth S^2_1}{2},\\
S^3_2=\frac{\overline{\goth S^1_2}\,\goth S^1_1
-\overline{\goth S^1_1}\,\goth S^1_2+\overline{\goth S^2_1}
\,\goth S^2_2-\overline{\goth S^2_2}\,\goth S^2_1}{2\,i},\\
S^3_3=\frac{\overline{\goth S^1_1}\,\goth S^1_1
+\overline{\goth S^2_2}\,\goth S^2_2-\overline{\goth S^2_1}
\,\goth S^2_1-\overline{\goth S^1_2}\,\goth S^1_2}{2}
\endgathered
\mytag{1.5}
\endgather
$$
(see \mycite{1}, \mycite{2}, and \mycite{3}, see also \mycite{4} and
\mycite{5} for more details). Here in \mythetag{1.2}, \mythetag{1.3}, 
\mythetag{1.4}, and \mythetag{1.5} by $\goth S^i_j$ we denote the
components of a $2\times 2$ complex matrix $\goth S\in\MatGrSL(2,
\Bbb C)$, while $S^i_j$ are the components of the matrix 
$S=\varphi(\goth S)$ produced from $\goth S$ by applying the 
homomorphism \mythetag{1.1}.\par
    Let $M$ be a {\it space-time} manifold, i\.\,e\. this is a 
$4$-dimensional orientable manifold equipped with a pseudo-Euclidean metric
$\bold g$ of the Minkowski-type signature $(+,-,-,-)$ and carrying a special
smooth geometric structure which is called a {\it polarization}. Once some
polarization is fixed, one can distinguish the {\it Future light cone\/}
from the {\it Past light cone\/} at each point $p\in M$ (see \mycite{6} for
more details). A moving frame $(U,\,\boldsymbol\Upsilon_0,\,\boldsymbol
\Upsilon_1,\,\boldsymbol\Upsilon_2,\,\boldsymbol\Upsilon_3)$ of the tangent 
bundle $TM$ is an ordered set of four smooth vector fields $\boldsymbol
\Upsilon_0$, $\boldsymbol\Upsilon_1$, $\boldsymbol\Upsilon_2$, $\boldsymbol
\Upsilon_3$ which are defined and $\Bbb R$-linearly independent at each
point $p$ of the open set $U\subset M$. This moving frame is called a {\it
positively polarized right orthonormal frame} if the following conditions 
are fulfilled:
\roster
\rosteritemwd=1pt
\item the value of the first vector filed $\boldsymbol\Upsilon_0$ at each
point $p\in U$ belongs to the interior of the Future light cone determined 
by the polarization of $M$;
\item it is a right frame in the sense of the orientation of $M$;
\item the metric tensor $\bold g$ is given by the standard
Minkowski matrix in this frame:
$$
\hskip -2em
g_{ij}=g(\boldsymbol\Upsilon_i,\boldsymbol\Upsilon_j)=\Vmatrix
\format \l&\quad\r&\quad\r&\quad\r\\
1 &0 &0 &0\\ 0 &-1 &0 &0\\ 0 &0 &-1 &0\\ 0 &0 &0 &-1\endVmatrix.
\mytag{1.6}
$$
\endroster
Assume that we have some other positively polarized right orthonormal frame
$(\tilde U,\,\tilde{\boldsymbol\Upsilon}_0,\,\tilde{\boldsymbol\Upsilon}_1,
\,\tilde{\boldsymbol\Upsilon}_2,\,\tilde{\boldsymbol\Upsilon}_3)$ such that
$U\cap\tilde U\neq\varnothing$. Then at each point $p\in U\cap\tilde U$ we
can write the following relationships for the frame vectors:
$$
\xalignat 2
&\hskip -2em
\tilde{\boldsymbol\Upsilon}_i=\sum^3_{j=0}S^j_i\,\boldsymbol\Upsilon_j,
&&\boldsymbol\Upsilon_i=\sum^3_{j=0}T^j_i\,\tilde{\boldsymbol\Upsilon}_j.
\mytag{1.7}
\endxalignat
$$
The relationships \mythetag{1.7} are called {\it transition formulas},
while the coefficients $S^j_i$ and $T^j_i$ in them are the components
of two mutually inverse transition matrices $S$ and $T$. Since both
frames $(U,\,\boldsymbol\Upsilon_0,\,\boldsymbol\Upsilon_1,\,\boldsymbol
\Upsilon_2,\,\boldsymbol\Upsilon_3)$ and $(\tilde U,\,\tilde{\boldsymbol
\Upsilon}_0,\,\tilde{\boldsymbol\Upsilon}_1,\,\tilde{\boldsymbol
\Upsilon}_2,\,\tilde{\boldsymbol\Upsilon}_3)$ are positively polarized
right orthonormal frames, the transition matrices $S$ and $T$ both are
orthochronous Lorentzian matrices with $\det S=1$ and $\det T=1$. Such
matrices form the {\it special orthochronous matrix Lorentz group}
$\MatGrSO^+(1,3,\Bbb R)$ (see \mycite{6} for more details).\par
    Let $SM$ be a two-dimensional smooth complex vector bundle over 
the space-time $M$ equipped with a non-vanishing skew-symmetric
spin-tensorial field $\bold d$. This spin-tensorial field $\bold d$ 
is called the {\it spin-metric tensor}. A moving frame $(U,\,\boldsymbol
\Psi_1,\,\boldsymbol\Psi_2)$ of $SM$ is an ordered set of two smooth 
section $\boldsymbol\Psi_1$ and $\boldsymbol\Psi_2$ of $SM$ over the
open set $U$ which are $\Bbb C$-linearly independent at each point 
$p\in U$. A moving frame $(U,\,\boldsymbol\Psi_1,\,\boldsymbol\Psi_2)$
is called an {\it orthonormal frame\/} if
$$
\hskip -2em
d_{ij}=d(\boldsymbol\Psi_i,\boldsymbol\Psi_j)
=\Vmatrix 0 & 1\\ 
\vspace{1ex} -1 & 0\endVmatrix,
\mytag{1.8}
$$
i\.\,e\. if the spin-metric tensor $\bold d$ is given by the skew-symmetric
matrix \mythetag{1.8} in this frame. Assume that we have two orthonormal
frames $(U,\,\boldsymbol\Psi_1,\,\boldsymbol\Psi_2)$ and $(\tilde U,\,
\tilde{\boldsymbol\Psi}_1,\,\tilde{\boldsymbol\Psi}_2)$ of the bundle 
$SM$ with overlapping domains $U\cap\tilde U\neq\varnothing$. Then at 
each point $p\in U\cap \tilde U$ we can write the following transition 
formulas:
$$
\xalignat 2
&\hskip -2em
\tilde{\boldsymbol\Psi}_i=\sum^2_{j=1}\goth S^j_i\,
\boldsymbol\Psi_j,
&&\boldsymbol\Psi_i=\sum^2_{j=1}\goth T^j_i\,
\tilde{\boldsymbol\Psi}_j.
\mytag{1.9}
\endxalignat
$$
Since both frames $(U,\,\boldsymbol\Psi_1,\,\boldsymbol\Psi_2)$ and
$(\tilde U,\,\tilde{\boldsymbol\Psi}_1,\,\tilde{\boldsymbol\Psi}_2)$
are orthonormal with respect to spin-metric tensor $\bold d$, both
transition matrices $\goth S$ and $\goth T=\goth S^{-1}$ with the
components $\goth S^j_i$ and $\goth T^j_i$ in \mythetag{1.9} belong
to the special linear group $\MatGrSL(2,\Bbb C)$. 
\mydefinition{1.1} A two-dimensional complex vector 
bundle $SM$ over the space-time manifold $M$ equipped with a nonzero
spin-metric $\bold d$ is called a {\it spinor bundle} if each
orthonormal frame $(U,\,\boldsymbol\Psi_1,\,\boldsymbol\Psi_2)$ of
$SM$ is associated with some positively polarized right orthonormal 
frame $(U,\,\boldsymbol\Upsilon_0,\,\boldsymbol\Upsilon_1,\,
\boldsymbol\Upsilon_2,\,\boldsymbol\Upsilon_3)$ of the tangent bundle
$TM$ such that for any two orthonormal frames $(U,\,\boldsymbol\Psi_1,
\,\boldsymbol\Psi_2)$ and $(\tilde U,\,\tilde{\boldsymbol\Psi}_1,\,
\tilde{\boldsymbol\Psi}_2)$ with overlapping domains $U\cap\tilde U\neq
\varnothing$ the associated tangent frames $(U,\,\boldsymbol\Upsilon_0,
\,\boldsymbol\Upsilon_1,\,\boldsymbol\Upsilon_2,\,\boldsymbol\Upsilon_3)$
and $(\tilde U,\,\tilde{\boldsymbol\Upsilon}_0,\,
\tilde{\boldsymbol\Upsilon}_1,\,\tilde{\boldsymbol\Upsilon}_2,\,
\tilde{\boldsymbol\Upsilon}_3)$ are related to each other by means of
the formulas \mythetag{1.7}, where the transition matrices $S$ and $T$
are obtained from the transition matrices $\goth S$ and $\goth T$ in
\mythetag{1.9} by applying the homomorphism \mythetag{1.1}, i\.\,e\. 
$S=\varphi(\goth S)$ and $T=\varphi(\goth T)$.
\enddefinition
\head
2. An algebraic background for Dirac spinors.
\endhead
     The group homomorphism \mythetag{1.1} is an algebraic background
for two-component spinors. They form a complex bundle over $M$ introduced
by the definition~\mythedefinition{1.1}. In order to construct an algebraic
background for Dirac spinors we need to extend the group homomorphism 
\mythetag{1.1} to bigger groups. For the group $\MatGrSO^+(1,3,\Bbb R)$
in \mythetag{1.1} we have the following natural enclosure:
$$
\hskip -2em
\MatGrSO^+(1,3,\Bbb R)\subset\MatGrO(1,3,\Bbb R).
\mytag{2.1}
$$
The complete matrix Lorentzian group $\MatGrO(1,3,\Bbb R)$ in
\mythetag{2.1} is generated by adding the following two matrices to
$\MatGrSO^+(1,3,\Bbb R)$:
$$
\xalignat 2
&\hskip -2em
P=\Vmatrix
\format \l&\quad\r&\quad\r&\quad\r\\
1 &0 &0 &0\\ 0 &-1 &0 &0\\ 0 &0 &-1 &0\\ 0 &0 &0 &-1\endVmatrix,
&&T=\Vmatrix
\format \l&\quad\r&\quad\r&\quad\r\\
-1 &0 &0 &0\\ 0 &1 &0 &0\\ 0 &0 &1 &0\\ 0 &0 &0 &1\endVmatrix.
\mytag{2.2}
\endxalignat
$$
The first matrix $P$ in \mythetag{2.2} is the {\it spatial inversion
matrix}. The second matrix $T$ is the {\it time inversion matrix}.
The matrices $P$ and $T$ in \mythetag{2.2} are commutating:
$$
\hskip -2em
P\cdot T=T\cdot P=-\bold 1.
\mytag{2.3}
$$
Apart from \mythetag{2.3}, we have the following relationships:
$$
\xalignat 2 
&\hskip -2em
P^2=\bold 1,&&T^2=\bold 1
\mytag{2.4}
\endxalignat
$$
Due to \mythetag{2.3} and \mythetag{2.4} each matrix $\hat S$ of 
the group $\MatGrO(1,3,\Bbb R)$ is represented as
$$
\hskip -2em
\hat S=S\text{\ \ \ or \ \ }
\hat S=P\cdot S\text{\ \ \ or \ \ }
\hat S=T\cdot S\text{\ \ \ or \ \ }
\hat S=-S,
\mytag{2.5}
$$
where $S\in\MatGrSO^+(1,3,\Bbb R)$. Using the representation \mythetag{2.5}, 
we can reduce the multiplication in $\MatGrO(1,3,\Bbb R)$ 
to the multiplication in $\MatGrSO^+(1,3,\Bbb R)$:
\bigskip
\hbox to \hsize{\hss
\vbox{\hsize 10cm
\offinterlineskip\settabs\+\indent
\vtrule
\hskip 1.1cm &\vttrule % Quantity
\hskip 2.3cm &\vtrule % Unit
\hskip 2.3cm &\vtrule % Unit
\hskip 2.3cm &\vtrule % Relation
\hskip 2.3cm &\vtrule % Relation
\cr\hrule 
\+\vtrule
\hss\hss &\vttrule
\hss $S_2$\hss &\vtrule
\hss $P\cdot S_2$\hss &\vtrule
\hss $T\cdot S_2$\hss &\vtrule
\hss $- S_2$\hss &\vtrule\cr
\hrule\hrule
\+\vtrule
\hss $S_1$\hss&\vttrule
\hss $S_1\cdot S_2$\hss&\vtrule
\hss $P\cdot(S_3\cdot S_2)$\hss&\vtrule
\hss $T\cdot(S_3\cdot S_2)$\hss&\vtrule
\hss $-(S_1\cdot S_2)$\hss&\vtrule\cr
\hrule
\+\vtrule
\hss $P\cdot S_1$\hss&\vttrule
\hss $P\cdot(S_1\cdot S_2)$\hss&\vtrule
\hss $S_3\cdot S_2$\hss&\vtrule
\hss $-(S_3\cdot S_2)$\hss&\vtrule
\hss $T\cdot(S_1\cdot S_2)$\hss&\vtrule\cr
\hrule
\+\vtrule
\hss $T\cdot S_1$\hss&\vttrule
\hss $T\cdot(S_1\cdot S_2)$\hss&\vtrule
\hss $-(S_3\cdot S_2)$\hss&\vtrule
\hss $S_3\cdot S_2$\hss&\vtrule
\hss $P\cdot(S_1\cdot S_2)$\hss&\vtrule\cr
\hrule
\+\vtrule
\hss $-S_1$\hss&\vttrule
\hss $-(S_1\cdot S_2)$\hss&\vtrule
\hss $T\cdot(S_3\cdot S_2)$\hss&\vtrule
\hss $P\cdot(S_3\cdot S_2)$\hss&\vtrule
\hss $S_1\cdot S_2$\hss&\vtrule\cr
\hrule
}\hss}
\bigskip\noindent
The matrix $S_3\in\MatGrSO^+(1,3,\Bbb R)$ in the above table is determined 
by the formula
$$
\hskip -2em
S_3=P\cdot S_1\cdot P.
\mytag{2.6}
$$
Since $P^2=1$, we have $P=P^{-1}$ and we can write \mythetag{2.6}
as $S_3=\psi(S_1)$, where $\psi\!:\,\MatGrSO^+(1,3,\Bbb R)\to
\MatGrSO^+(1,3,\Bbb R)$ is the group homomorphism given by the
formula
$$
\hskip -2em
S\mapsto S'=\psi(S)=P\cdot S\cdot P^{-1}.
\mytag{2.7}
$$
Assume that the matrix $S$ is obtained by means of the group homomorphism
\mythetag{1.1}, i\.\,e\. assume that $S=\varphi(\goth S)$, where $\goth S
\in\MatGrSL(2,\Bbb C)$. Then $S\in\MatGrSO^+(1,3,\Bbb R)$, and
$\psi(S)\in\MatGrSO^+(1,3,\Bbb R)$, therefore $\psi(S)=\varphi(\goth S')$
for some matrix $\goth S\in\MatGrSL(2,\Bbb C)$ since $\varphi$ is a surjective 
homomorphism. The matrix $\goth S'$ is determined uniquely
up to the sign. By means of direct calculations we derive
$$
\hskip -2em
\goth S'=\pm\,(\goth S^{-1})^{\sssize\dagger},
\mytag{2.8}
$$
i\.\,e\. $\goth S'$ is the Hermitian conjugate matrix for $\goth S^{-1}$.
Choosing the plus sign in \mythetag{2.8}, we obtain one more group
homomorphism
$$
\hskip -2em
\goth S\mapsto\goth S'=\psi'(\goth S)=(\goth S^{-1})^{\sssize\dagger}.
\mytag{2.9}
$$
The homomorphisms \mythetag{1.1}, \mythetag{2.7}, and \mythetag{2.9}
compose the commutative diagram
$$
\CD
\goth S @>\varphi>>S\\
@V\psi'VV @VV\psi V\\
\goth S'@>\varphi>>S'
\endCD
$$\par
     Now let's remember the construction of the homomorphism \mythetag{1.1},
see \mycite{1}, \mycite{2}, \mycite{3}, or \mycite{4}.
It is constructed on the base of the equality 
$$
\hskip -2em
\goth S\cdot\boldsymbol\sigma_m\cdot\goth S^{\sssize\dagger}=\sum^3_{k=0}
S^k_m\,\boldsymbol\sigma_k,
\mytag{2.10}
$$
where $\boldsymbol\sigma_0$ is the unit matrix, while 
$\boldsymbol\sigma_1,\,\boldsymbol\sigma_2,\,\boldsymbol\sigma_3$ are
the well-known Pauli matrices:
$$
\xalignat 2
&\hskip -2em
\boldsymbol\sigma_0=\Vmatrix 1 & 0\\0 & 1\endVmatrix,
&&\boldsymbol\sigma_2=\Vmatrix 0 & -i\\i & 0\endVmatrix,\\
\vspace{-1.4ex}
&&&\mytag{2.11}\\
\vspace{-1.4ex}
&\hskip -2em
\boldsymbol\sigma_1=\Vmatrix 0 & 1\\1 & 0\endVmatrix,
&&\boldsymbol\sigma_3=\Vmatrix 1 & 0\\0 & -1\endVmatrix.
\endxalignat
$$
\mylemma{2.1} For $\goth S\in\MatGrSL(2,\Bbb C)$ the relationship
\mythetag{2.10} can be transformed to
$$
\hskip -2em
(\goth S^{-1})^{\sssize\dagger}\cdot\tilde{\boldsymbol\sigma}_m\cdot
\goth S^{-1}=\sum^3_{k=0}S^k_m\,\tilde{\boldsymbol\sigma}_k,
\mytag{2.12}
$$
where $\tilde{\boldsymbol\sigma}_m=\varepsilon_m\,\boldsymbol\sigma_m^{-1}$
and $\varepsilon_m=\det(\boldsymbol\sigma_m)$, i\.\,e\. they are given
by the formulas
$$
\xalignat 2
&\hskip -2em
\tilde{\boldsymbol\sigma}_0=\Vmatrix 1 & 0\\0 & 1\endVmatrix,
&&\tilde{\boldsymbol\sigma}_2=\Vmatrix 0 & i\\-i & 0\endVmatrix,\\
\vspace{-1.4ex}
&&&\mytag{2.13}\\
\vspace{-1.4ex}
&\hskip -2em
\tilde{\boldsymbol\sigma}_1=\Vmatrix 0 & -1\\-1 & 0\endVmatrix,
&&\tilde{\boldsymbol\sigma}_3=\Vmatrix -1 & 0\\0 & 1\endVmatrix.
\endxalignat
$$
\endproclaim
\demo{Proof} Note that for a $2\times 2$ matrix $A$ with the unit
determinant $\det A=1$ we have
$$
\hskip -2em
A^{-1}=\Vmatrix a^1_1 & a^1_2\\\vspace{1ex}a^2_1 & a^2_2\endVmatrix^{-1}
=\Vmatrix a^2_2 & -a^1_2\\\vspace{1ex}-a^2_1 & a^1_1\endVmatrix.
\mytag{2.14}
$$
Relying on \mythetag{2.14}, we define the map taking a $2\times 2$ matrix 
to another $2\times 2$ matrix:
$$
\hskip -2em
A\mapsto L(A)=L\!\left(\,\Vmatrix a^1_1 & a^1_2\\\vspace{1ex}a^2_1 &
a^2_2\endVmatrix\,\right)=\Vmatrix a^2_2 & -a^1_2\\\vspace{1ex}-a^2_1
& a^1_1\endVmatrix.
\mytag{2.15}
$$
If $\det A=\pm\,1$, then for this $2\times 2$ matrix $A$ we have
$$
\hskip -2em
A^{-1}=\det(A)\cdot L(A).
\mytag{2.16}
$$
Note that $\varepsilon_m=\det(\boldsymbol\sigma_m)=\pm\,1$ for all matrices
\mythetag{2.11}. The same is true for the matrix in the left hand side
of \mythetag{2.10} since $\goth S\in\MatGrSL(2,\Bbb C)$ and $\det\goth S=1$.
Therefore, applying \mythetag{2.16} to \mythetag{2.10}, we derive
$$
\hskip -2em
(\goth S^{-1})^{\sssize\dagger}\cdot{\boldsymbol\sigma}_m^{-1}\cdot
\goth S^{-1}=\varepsilon_m\ L\!\left(\,\shave{\sum^3_{k=0}}S^k_m\,
\boldsymbol\sigma_k\right).
\mytag{2.17}
$$
It is easy to see that the map \mythetag{2.15} is a linear map. For this
reason we can transform the above equality \mythetag{2.17} in the following
way:
$$
\hskip -2em
(\goth S^{-1})^{\sssize\dagger}\cdot\boldsymbol\sigma_m^{-1}\cdot
\goth S^{-1}=\varepsilon_m\ \sum^3_{k=0}S^k_m\,L(\boldsymbol\sigma_k)=
\varepsilon_m\ \sum^3_{k=0}S^k_m\ \varepsilon_k\ \sigma_k^{-1}.
\mytag{2.18}
$$
Since $\varepsilon_m=\pm\,1$ and $\varepsilon_k=\pm\,1$, passing from
\mythetag{2.11} to the matrices \mythetag{2.13} in \mythetag{2.18}, we
see that it coincides with the required equality \mythetag{2.12}.
\qed\enddemo
     The next step is to extend the group $\MatGrSL(2,\Bbb C)$ using the
homomorphism \mythetag{2.9} for this purpose. It is easy to see that the
mapping 
$$
\hskip -2em
\goth S\mapsto\Psi'(\goth S)=\hat{\goth S}
=\Vmatrix\format\ \c\ &\,\c\ &\c\\
\goth S &\vrule height 8pt 
depth 4pt&0\\
\vspace{-8pt}
\kern -14pt\vbox{\hrule width 50pt}\kern -50pt&&\\
\vspace{-1pt}
0 &\vrule height 12pt depth 3pt
&(\goth S^{-1})^{\sssize\dagger}
\endVmatrix
\mytag{2.19}
$$
is an embedding of the group $\MatGrSL(2,\Bbb C)$ into the general linear
group $\MatGrGL(4,\Bbb C)$, i\.\,e\. it is an exact representation of the
group $\MatGrSL(2,\Bbb C)$ by means of complex $4\times 4$ matrices. In
addition to \mythetag{2.19}, we consider the following $4\times 4$ matrices:
$$
\hskip -2em
\gamma_m=\Vmatrix\format\c&\,\c\ &\c\\
0 &\vrule height 8pt 
depth 4pt&\boldsymbol\sigma_m\\
\vspace{-8pt}
\kern -14pt\vbox{\hrule width 34pt}\kern -34pt&&\\
\vspace{-1pt}
\tilde{\boldsymbol\sigma}_m&\vrule height 12pt depth 3pt
&0\endVmatrix\text{, \ \ }m=0,1,2,3.
\mytag{2.20}
$$
Then, using the matrices \mythetag{2.19} and \mythetag{2.20}, we combine
them as follows:
$$
\hskip -2em
\hat{\goth S}\cdot\gamma_m\cdot\hat{\goth S}^{-1}=
\Vmatrix\format\c&\,\c\ &\c\\
0 &\vrule height 8pt 
depth 4pt&\goth S\cdot\boldsymbol\sigma_m\cdot\goth S^{\sssize\dagger}\\
\vspace{-8pt}
\kern -78pt\vbox{\hrule width 133pt}\kern -133pt&&\\
\vspace{-1pt}
(\goth S^{-1})^{\sssize\dagger}\cdot
\tilde{\boldsymbol\sigma}_m\cdot\goth S^{-1}
&\vrule height 12pt depth 3pt
&0\endVmatrix.
\mytag{2.21}
$$
Applying \mythetag{2.10} and \mythetag{2.18} to \mythetag{2.21}, we can
transform this equality to
$$
\hskip -2em
\hat{\goth S}\cdot\gamma_m\cdot\hat{\goth S}^{-1}=
\sum^3_{k=0}S^k_m\,\gamma_k.
\mytag{2.22}
$$\par
    The matrices $\gamma_m$ in \mythetag{2.20} are known as Dirac matrices.
They obey the following anticommutation relationship (see \mycite{2}):
$$
\hskip -2em
\{\gamma_i,\,\gamma_j\}=2\,g_{ij}\ \bold 1.
\mytag{2.23}
$$
Here $\{\gamma_i,\,\gamma_j\}=\gamma_i\cdot\gamma_j-\gamma_j\cdot\gamma_i$
is the matrix anticommutator, while $\bold 1$ in \mythetag{2.23} is the unit
$4\times 4$ matrix and $g_{ij}$ are the components of the matrix 
\mythetag{1.6}, i\.\,e\. they are numbers. Due to \mythetag{2.23} the
following combinations of $\gamma$-matrices are independent:
$$
\allowdisplaybreaks
\xalignat 2
\hskip -2em
\bold 1&=\Vmatrix 1&0&0&0\\ 0&1&0&0\\ 0&0&1&0\\ 0&0&0&1\endVmatrix,
&\kern -40pt\gamma_0\cdot\gamma_1\cdot\gamma_2\cdot\gamma_3&=\Vmatrix
i&0&0&0\\0&i&0&0\\0&0&-i&0\\0&0&0&-i\endVmatrix,\qquad
\mytag{2.24}\\
\hskip -2em
\gamma_0&=\Vmatrix 0&0&1&0\\ 0&0&0&1\\ 1&0&0&0\\ 0&1&0&0\endVmatrix,
&\kern -35pt\gamma_1\cdot\gamma_2\cdot\gamma_3&=\Vmatrix 0&0&-i&0\\
0&0&0&-i\\i&0&0&0\\ 0&i&0&0\endVmatrix,
\mytag{2.25}\\
\hskip -2em
\gamma_1&=\Vmatrix 0&0&0&1\\ 0&0&1&0\\ 0&-1&0&0\\ -1&0&0&0\endVmatrix,
&\kern -35pt\gamma_0\cdot\gamma_2\cdot\gamma_3&=\Vmatrix 0&0&0&-i\\
0&0&-i&0\\0&-i&0&0\\ -i&0&0&0\endVmatrix,
\mytag{2.26}\\
\hskip -2em
\gamma_2&=\Vmatrix 0&0&0&-i\\ 0&0&i&0\\ 0&i&0&0\\ -i&0&0&0\endVmatrix,
&\kern -35pt\gamma_0\cdot\gamma_1\cdot\gamma_3&=\Vmatrix 0&0&0&1\\
0&0&-1&0\\0&1&0&0\\ -1&0&0&0\endVmatrix,
\mytag{2.27}\\
\hskip -2em
\gamma_3&=\Vmatrix 0&0&1&0\\ 0&0&0&-1\\ -1&0&0&0\\ 0&1&0&0\endVmatrix,
&\kern -35pt\gamma_0\cdot\gamma_1\cdot\gamma_2&=\Vmatrix 0&0&-i&0\\
0&0&0&i\\-i&0&0&0\\0&i&0&0\endVmatrix,
\mytag{2.28}\\
\hskip -4em
\gamma_0\cdot\gamma_1&=\Vmatrix 0&-1&0&0\\ -1&0&0&0\\ 0&0&0&1\\
0&0&1&0\endVmatrix,
&\kern -25pt\gamma_2\cdot\gamma_3&=\Vmatrix 0&-i&0&0\\ -i&0&0&0\\ 
0&0&0&-i\\0&0&-i&0\endVmatrix,
\mytag{2.29}\\
\hskip -4em
\gamma_0\cdot\gamma_2&=\Vmatrix 0&i&0&0\\ -i&0&0&0\\ 0&0&0&-i\\
0&0&i&0\endVmatrix,
&\kern -35pt\gamma_1\cdot\gamma_3&=\Vmatrix 0&1&0&0\\-1&0&0&0\\ 
0&0&0&1\\0&0&-1&0\endVmatrix,
\mytag{2.30}\\
\hskip -4em
\gamma_0\cdot\gamma_3&=\Vmatrix -1&0&0&0\\ 0&1&0&0\\ 0&0&1&0\\
0&0&0&-1\endVmatrix,
&\kern -35pt\gamma_1\cdot\gamma_2&=\Vmatrix -i&0&0&0\\ 0&i&0&0\\ 
0&0&-i&0\\0&0&0&i\endVmatrix,
\mytag{2.31}\\
\endxalignat
$$
All other products of Dirac matrices are expressed as linear combinations
of these matrices. Moreover, the matrices \mythetag{2.24}, \mythetag{2.25},
\mythetag{2.26}, \mythetag{2.27}, \mythetag{2.28}, \mythetag{2.29},
\mythetag{2.30}, \mythetag{2.31} are linearly independent over the field
of complex numbers $\Bbb C$. They form a basis in the linear space of
$4\times 4$ complex matrices. This fact is well-known, it is mentioned
in \mycite{2}.\par
     In order to extend the group $\MatGrSL(2,\Bbb C)$ represented as 
a subgroup $G\subset\MatGrGL(4,\Bbb C)$ by means of the embedding
\mythetag{2.19} we modify the relationship \mythetag{2.22} as follows.
We replace $S^k_m$ in \mythetag{2.22} by the components of the spatial
inversion matrix $P$ from \mythetag{2.2}. Then we replace $\hat{\goth S}$
by some $4\times 4$ matrix $\hat P$ which is yet unknown:
$$
\hskip -2em
\hat P\cdot\gamma_m\cdot\hat P^{-1}=
\sum^3_{k=0}P^k_m\,\gamma_k.
\mytag{2.32}
$$
Our next goal is to solve \mythetag{2.32} with respect to the unknown
matrix $\hat P$. Note that it can be written as a system of four
very simple matrix equations:
$$
\xalignat 2
&\hat P\cdot\gamma_0=\gamma_0\cdot\hat P,
&&\hat P\cdot\gamma_1=-\gamma_1\cdot\hat P,\\
&\hat P\cdot\gamma_2=-\gamma_2\cdot\hat P,
&&\hat P\cdot\gamma_3=-\gamma_3\cdot\hat P.
\endxalignat
$$
This is the system of $64$ linear homogeneous equations with respect to
16 components of the matrix $\hat P$. Its general solution is given by
the formula
$$
\hskip -2em
\hat P=C\,\gamma_0,
\mytag{2.33}
$$
where $C$ is an arbitrary complex number.\par
     In a similar way, taking the components of time inversion matrix 
$T$ from \mythetag{2.2}, on the base of \mythetag{2.22} we can write 
the equation
$$
\hskip -2em
\hat T\cdot\gamma_m\cdot\hat T^{-1}
=\sum^3_{k=0}T^k_m\,\gamma_k
\mytag{2.34}
$$
for the unknown matrix $\hat T$. Like the equation \mythetag{2.32}, the 
equation \mythetag{2.34} reduces to a system of four matrix equations:
$$
\xalignat 2
&\hat T\cdot\gamma_0=-\gamma_0\cdot\hat T,
&&\hat T\cdot\gamma_1=\gamma_1\cdot\hat T,\\
&\hat T\cdot\gamma_2=\gamma_2\cdot\hat T,
&&\hat T\cdot\gamma_3=\gamma_3\cdot\hat T.
\endxalignat
$$
The general solution of this system of equations is given by the formula
$$
\hskip -2em
T=C\,\gamma_1\cdot\gamma_2\cdot\gamma_3,
\mytag{2.35}
$$
where $C$ again is an arbitrary complex number. In order to fix the
constants in \mythetag{2.33} and \mythetag{2.35} we apply the following
restrictions\footnote{\, In some cases other normalization conditions
for $\hat P$ and $\hat T$ are used, see \mycite{7} and \mycite{8}.} to
$\hat P$ and $\hat T$:
\adjustfootnotemark{-1}
$$
\xalignat 2
&\hskip -2em
\hat P^2=\bold 1, &&\hat T^2=\bold 1.
\mytag{2.36}
\endxalignat
$$
From \mythetag{2.36}, \mythetag{2.33}, and \mythetag{2.35} we derive
$$
\xalignat 2
&\hskip -2em
\hat P=\pm\,\gamma_0,
&&\hat T=\pm\,\gamma_1\cdot\gamma_2\cdot\gamma_3.
\mytag{2.37}
\endxalignat
$$
Moreover, from \mythetag{2.36} and \mythetag{2.37} we derive
$$
\xalignat 2
&\hskip -2em
\{\hat P,\,\hat T\}=0,
&&(\hat P\cdot\hat T)^2=-\bold 1.
\mytag{2.38}
\endxalignat
$$
The first equality \mythetag{2.38} means that $\hat P$ and $\hat T$
are anticommutative with respect to each other. The relationships 
\mythetag{2.38} are valid for any choice of sign in \mythetag{2.37}.
This fact can be strengthened in the following way.
\mylemma{2.2}For any choice of signs the matrices \mythetag{2.37} and 
the matrices \mythetag{2.19}, where $\goth S\in\MatGrSL(2,\Bbb C)$, 
generate the same subgroup $G\subset\MatGrGL(4,\Bbb C)$ being
a discrete extension of the \pagebreak group $\MatGrSL(2,\Bbb C)$.
\endproclaim
    The group $G$ in the lemma~\mythelemma{2.2} is isomorphic to the
spinor group $\GrPin(1,3,\Bbb R)$. If we identify $G$ with 
$\GrPin(1,3,\Bbb R)$ according to this isomorphism, then the subgroup
\footnote{\,The definitions of the groups $\GrSpin(1,3,\Bbb R)$ and
$\GrPin(1,3,\Bbb R)$ can be found in \mycite{9}.} $\GrSpin(1,3,\Bbb R)
\subset\GrPin(1,3,\Bbb R)$ is identified with the $4$-dimensional
presentation of the group $\MatGrSL(2,\Bbb C)$ given by the matrices 
\mythetag{2.19}. Thus we have reached the goal stated in the very 
beginning of this section. By introducing the matrices \mythetag{2.19} 
and \mythetag{2.37} we have constructed the group homomorphism 
\adjustfootnotemark{-1}
$$
\hskip -2em
\Phi\!:\,\GrPin(1,3,\Bbb R)\to\MatGrO(1,3,\Bbb R).
\mytag{2.39}
$$
This homomorphism \mythetag{2.39} extends the initial homomorphism 
\mythetag{1.1} in the sense of the following commutative diagram:
$$
\hskip -2em
\CD
\MatGrSL(2,\Bbb C) @>\varphi>>\MatGrSO(1,3,\Bbb R)\\
@V\Psi'VV @VVV\\
\kern -62pt\GrPin(1,3,\Bbb R)\cong G@>\Phi>>\MatGrO(1,3,\Bbb R).
\endCD
\mytag{2.40}
$$
Both vertical arrows in the diagram \mythetag{2.40} are embeddings. The
homomorphism \mythetag{2.39} is described by the formula
$$
\hskip -2em
\hat{\goth S}\cdot\gamma_m\cdot\hat{\goth S}^{-1}=
\sum^3_{k=0}S^k_m\,\gamma_k.
\mytag{2.41}
$$
This formula coincides with \mythetag{2.22}, however, now $\hat{\goth S}$
is not necessarily given by the formula \mythetag{2.19}. It is an arbitrary
matrix from the group $G\cong\GrPin(1,3,\Bbb R)$. In particular, we can 
take $\hat{\goth S}=\hat P$ or $\hat{\goth S}=\hat T$. Then \mythetag{2.41}
reduces to \mythetag{2.32} or to \mythetag{2.34} respectively. By $S^k_m$
now in the formula \mythetag{2.41} we denote the components of the
Lorentzian matrix $S=\Phi(\hat{\goth S})\in\MatGrO(1,3,\Bbb R)$.\par
     Like \mythetag{1.1}, the homomorphism \mythetag{2.39} is a srjective
mapping. Its kernel is discrete, it is composed by the following two
matrices:
$$
\xalignat 2
&\hskip -2em
\bold 1=\Vmatrix 1&0&0&0\\ 0&1&0&0\\ 0&0&1&0\\ 0&0&0&1\endVmatrix,
&&-\bold 1=\Vmatrix -1&0&0&0\\ 0&-1&0&0\\ 0&0&-1&0\\ 0&0&0&-1\endVmatrix.
\quad
\mytag{2.42}
\endxalignat
$$
Due to \mythetag{2.42} for any matrix $S\in\MatGrO(1,3,\Bbb R)$ its
preimage $\hat{\goth S}\in G\cong\GrPin(1,3,\Bbb R)$ is determined 
uniquely up to the sign.
\head
3. Dirac spinors.
\endhead
    Let $M$ be a space-time manifold and let $SM$ be a spinor bundle
over $M$ introduced by the definition~\mythedefinition{1.1}. By 
$S^{\sssize\dagger}\!M$ we denote the Hermitian conjugate bundle for
$SM$. Taking both $SM$ and $S^{\sssize\dagger}\!M$, we construct their
direct sum
$$
\hskip -2em
DM=SM\oplus S^{\sssize\dagger}\!M.
\mytag{3.1}
$$
The direct sum \mythetag{3.1} is called the {\it Dirac bundle} associated
with the spinor bundle $SM$. This is a four-dimensional complex bundle
over $M$. The bundles $SM$ and $S^{\sssize\dagger}\!M$, when treated as
the constituents of $DM$, are called {\it chiral bundles}. Local and
global smooth sections of the Dirac bundle $DM$ are called {\it spinor 
fields} or more precisely {\it spinor fields of Dirac spinors}.\par
     According to the definition~\mythedefinition{1.1}, the chiral bundle
$SM$ in \mythetag{3.1} is equipped with the spin-metric $\bold d$. This
metric induces dual metric $\bold d$ in $S^*\!M$. Then by means of the
semilinear isomorphism of complex conjugation 
$$
\hskip -2em
\tau\!:\,S^*\!M\to S^{\sssize\dagger}\!M
\mytag{3.2}
$$
it is transferred to the Hermitian conjugate bundle $S^{\sssize\dagger}
\!M$ (see more details in \mycite{4}). Having spin-metrics in $SM$ and 
in $S^{\sssize\dagger}\!M$, we can define a spin-metric in $DM$. Indeed,
let $\bold X\in D_p(M)$ and $\bold Y\in D_p(M)$ at some point $p\in M$.
Then due to \mythetag{3.1} we have
$$
\align
&\hskip -2em
\bold X=\bold X_1+\bold X_2\text{, \ where \ }\bold X_1\in S_p(M)
\text{\ \ and \ }\bold X_2\in S^{\sssize\dagger}_p(M),\\
\vspace{-1.5ex}
\mytag{3.3}\\
\vspace{-1.5ex}
&\hskip -2em
\bold Y=\bold Y_1+\bold Y_2\text{, \ where \ }\bold Y_1\in S_p(M)
\text{\ \ and \ }\bold Y_2\in S^{\sssize\dagger}_p(M).
\endalign
$$
Using the expansions \mythetag{3.3}, by definition we set
$$
\hskip -2em
d(\bold X,\bold Y)=d(\bold X_1,\bold Y_1)+d(\bold X_2,\bold Y_2).
\mytag{3.4}
$$
Thus, the skew-symmetric metric \mythetag{3.4} in $DM$ is introduced 
as the sum of metrics in $SM$ and $S^{\sssize\dagger}\!M$ due to the
expansion \mythetag{3.1}.\par
     For each vector $\bold X\in D_p(M)$ we have the expansion
\mythetag{3.3} determined by the expansion \mythetag{3.1}. The 
operator $\bold H$ then is defined by the formula
$$
\hskip -2em
\bold H(\bold X)=\bold X_1-\bold X_2.
\mytag{3.5}
$$
This formula means that $\bold H$ in each fiber $D_p(M)$ is defined as 
a linear operator with two eigenvalues $\lambda_1=1$ and $\lambda_2=-1$. 
The eigenspace for $\lambda_1$ coincides with $S_p(M)$ and the 
eigenspace for $\lambda_2$ coincides with $S^{\sssize\dagger}_p(M)$. The
operator field $\bold H$ introduced by the formula \mythetag{3.5} is
called the {\it chirality operator}.\par
     The inverse map for \mythetag{3.2} is denoted by the same symbol
$\tau$. It is also called the semilinear isomorphism of complex conjugation
(see \mycite{4}):
$$
\hskip -2em
\tau\!:\,S^{\sssize\dagger}\!M\to S^*\!M.
\mytag{3.6}
$$
Applying \mythetag{3.6} to $\bold X_2$ and $\bold Y_2$ in \mythetag{3.3},
we get two chiral cospinors in $S^*_p(M)$:
$$
\xalignat 2
&\hskip -2em
\bold x_2=\tau(\bold X_2),
&&\bold y_2=\tau(\bold Y_2).
\mytag{3.7}
\endxalignat
$$
Since $\bold x_2\in S^*_p(M)$ and $\bold x_2\in S^*_p(M)$, they can be
paired with $\bold X_1$ and $\bold Y_1$. Hence, we can define the following
pairing for the spinors $\bold X$ and $\bold Y$:
$$
\hskip -2em
D(\bold X,\bold Y)=(\bold x_2,\bold Y_1)+\overline{(\bold y_2,\bold X_1)}
\mytag{3.8}
$$
The Hermitian form $\bold D$ defined by means of the formulas \mythetag{3.7}
and \mythetag{3.8} is called the {\it Dirac form} or the {\it Hermitian
spin-metric}. Note that the Hermitian spin-metric $\bold D$ is not positive,
its signature \pagebreak is $(+,+,-,-)$.\par
     The spin-metric $\bold d$, the chirality operator $\bold H$ and 
the Hermitian spin-metric $\bold D$ are basic geometric structures
associated with Dirac spinors. Some other structures will be considered
below a little bit later.\par
\head
4. Spin-tensors.
\endhead
     The definition of spin-tensors in the case of Dirac spinors is quite
standard. We introduce them following the scheme of the paper \mycite{4}.
Let $T_p(M)$ and $D_p(M)$ be the fibers of the tangent bundle $TM$ and 
the Dirac bundle $DM$ at some point $p\in M$. Denote by $T^*_p(M)$ and 
$D^*_p(M)$ the dual spaces for $T_p(M)$ and $D_p(M)$, then produce 
from $T_p(M)$ and $T^*_p(M)$ the complex spaces $\Bbb CT_p(M)$ and 
$\Bbb C T^*_p(M)$ by means of standard complexification procedure:
$$
\xalignat 2
&\hskip -2em
\Bbb CT_p(M)=\Bbb C\otimes T_p(M),
&&\Bbb CT^*_p(M)=\Bbb C\otimes T^*_p(M).
\mytag{4.1}
\endxalignat
$$
The complex spaces \mythetag{4.1} are obviously dual to each other.
In addition to $D_p(M)$ and $D^*_p(M)$ we introduce the Hermitian conjugate
spaces 
$$
\xalignat 2
&\hskip -2em
D^{\sssize\dagger}_p(M),
&&D^{*\sssize\dagger}_p(M)=D^{{\sssize\dagger}*}_p(M).
\mytag{4.2}
\endxalignat
$$
Then, using \mythetag{4.1} and \mythetag{4.2}, we define the following
tensor products:
$$
%\allowdisplaybreaks
\gather
\hskip -2em
\Bbb CT^m_n(p,M)=\overbrace{\Bbb CT_p(M)\otimes\ldots\otimes
\Bbb CT_p(M)}^{\text{$m$ times}}\otimes\underbrace{\Bbb CT^*_p(M)
\otimes\ldots\otimes \Bbb CT^*_p(M)}_{\text{$n$ times}},
\mytag{4.3}\\
\hskip -2em
D^\alpha_\beta(p,M)=\overbrace{D_p(M)\otimes\ldots\otimes 
D_p(M)}^{\text{$\alpha$ times}}\otimes
\underbrace{D^*_p(M)\otimes\ldots\otimes 
D^*_p(M)}_{\text{$\beta$ times}},
\mytag{4.4}\\
%\displaybreak
\hskip -2em
\bar D^\nu_\gamma(p,M)=\overbrace{D^{{\sssize\dagger}*}_p(M)\otimes
\ldots\otimes D^{{\sssize\dagger}*}_p(M)}^{\text{$\nu$ times}}
\otimes\underbrace{D^{\sssize\dagger}_p(M)\otimes\ldots\otimes 
D^{\sssize\dagger}_p(M)}_{\text{$\gamma$ times}}.
\mytag{4.5}
\endgather
$$     
Combining \mythetag{4.3}, \mythetag{4.4}, and \mythetag{4.5}, we define
one more tensor product
$$
\hskip -2em
D^\alpha_\beta\bar D^\nu_\gamma T^r_s(p,M)=D^\alpha_\beta(p,M)
\otimes\bar D^\nu_\gamma(p,M)\otimes\Bbb CT^m_n(p,M).
\mytag{4.6}
$$
Elements of the space \mythetag{4.6} are called {\it Dirac spin-tensors\/} of
the type $(\alpha,\beta|\nu,\gamma|m,n)$ at the point $p\in M$. The spaces
\mythetag{4.6} with $p$ running over the space-time manifold $M$ are
naturally glued into a bundle. This bundle is called the {\it spin-tensorial
bundle\/} of the type $(\alpha,\beta|\nu,\gamma|m,n)$, its local and global
smooth sections are called {\it spin-tensorial fields\/} of the type
$(\alpha,\beta|\nu,\gamma|m,n)$.\par
     The complex conjugation isomorphism $\tau$ for Dirac spin-tensors
is introduced in the same way as in the case of chiral spin-tensors.
The tangent space $T_p(M)$ and the cotangent space $T^*_p(M)$ are real
spaces. Therefore, here we have
$$
\hskip -2em
\aligned
&\tau(\bold X)=\bold X\text{\ \ for all \ }\bold X\in T_p(M),\\
&\tau(\lambda\ \bold X)=\overline{\lambda}\ \bold X\text{\ \ for any \ }
\lambda\in\Bbb C.
\endaligned
\mytag{4.7}
$$
Similarly, in the case of the cotangent space $T^*_p(M)$ we have
$$
\hskip -2em
\aligned
&\tau(\bold u)=\bold u\text{\ \ for all \ }\bold u\in T^*_p(M),\\
&\tau(\lambda\ \bold u)=\overline{\lambda}\ \bold u\text{\ \ for any \ }
\lambda\in\Bbb C.
\endaligned
\mytag{4.8}
$$
The formulas \mythetag{4.7} and \mythetag{4.8} define $\tau$ as two
semilinear mappings
$$
\xalignat 2
&\hskip -2em
\tau\!:\,\Bbb CT_p(M)\to\Bbb CT_p(M),
&&\tau\!:\,\Bbb CT^*_p(M)\to\Bbb CT^*_p(M)
\mytag{4.9}
\endxalignat
$$
such that $\tau^2=\tau\compos\tau=\idop$. The mappings \mythetag{4.9} are
easily extended to the tensor product \mythetag{4.3}. As a result we have
the semilinear mapping 
$$
\hskip -2em
\tau\!:\,\Bbb CT^m_n(p,M)\to\Bbb CT^m_n(p,M)
\mytag{4.10}
$$
with the same property $\tau^2=\tau\compos\tau=\idop$. Apart from the
mappings \mythetag{4.9}, we have the following canonical semilinear
mappings mutually inverse in each pair:
$$
\xalignat 2
&\hskip -2em
\CD
@>\tau>>\\
\vspace{-4ex}
D_p(M)@.D^{*\sssize\dagger}_p(M),\\
\vspace{-4.2ex}
@<<\tau< 
\endCD
&&\CD
@>\tau>>\\
\vspace{-4ex}
D^*_p(M)@.D^{\sssize\dagger}_p(M).\\
\vspace{-4.2ex}
@<<\tau< 
\endCD\quad
\mytag{4.11}
\endxalignat
$$
They are defined according to the recipe of the section~3 in \mycite{4}. 
All of the mappings \mythetag{4.11} are denoted by the same symbol 
$\tau$ so that we formally preserve the property $\tau^2=\tau\compos
\tau=\idop$. They are easily extended to the tensor products 
\mythetag{4.4} and \mythetag{4.5}:
$$
\hskip -2em
\CD
@>\tau>>\\
\vspace{-4ex}
D^\alpha_\beta(p,M)@.\bar D^\alpha_\beta(p,M),\\
\vspace{-4.2ex}
@<<\tau< 
\endCD
\mytag{4.12}
$$
Note that the second pair of the mappings \mythetag{4.11} are similar
to \mythetag{3.2} and \mythetag{3.6}, though here they are defined
independently. Now, combining the mappings \mythetag{4.10} and
\mythetag{4.12}, we extend $\tau$ to the tensor product \mythetag{4.6}:
$$
\hskip -2em
\CD
@>\tau>>\\
\vspace{-4ex}
D^\alpha_\beta\bar D^\nu_\gamma T^r_s(p,M)@.
D^\nu_\gamma\bar D^\alpha_\beta T^r_s(p,M).\\
\vspace{-4.2ex}
@<<\tau< 
\endCD
\mytag{4.13}
$$
The mappings \mythetag{4.13} are inverse to each other so that the
property $\tau^2=\tau\compos\tau=\idop$ for them is again formally
preserved.\par
     Let $(U,\,\boldsymbol\Psi_1,\,\boldsymbol\Psi_2,\,\boldsymbol\Psi_3,
\,\boldsymbol\Psi_4)$ and $(\tilde U,\,\tilde{\boldsymbol\Psi}_1,\,
\tilde{\boldsymbol\Psi}_2,\,\tilde{\boldsymbol\Psi}_3,\,\tilde{\boldsymbol
\Psi}_4)$ be two frames of the Dirac bundle $DM$ with overlapping domains:
$U\cap\tilde U\neq\varnothing$. Let $(U,\,\boldsymbol\Upsilon_0,\,
\boldsymbol\Upsilon_1,\,\boldsymbol\Upsilon_2,\,\boldsymbol\Upsilon_3)$
and $(\tilde U,\,\tilde{\boldsymbol\Upsilon}_0,\,\tilde{\boldsymbol
\Upsilon}_1,\,\tilde{\boldsymbol\Upsilon}_2,\,\tilde{\boldsymbol
\Upsilon}_3)$ be two frames of the tangent bundle $TM$ with the same 
domains $U$ and $\tilde U$. Note that here we do not require $(U,\,
\boldsymbol\Upsilon_0,\,\boldsymbol\Upsilon_1,\,\boldsymbol\Upsilon_2,
\,\boldsymbol\Upsilon_3)$ and $(\tilde U,\,\tilde{\boldsymbol\Upsilon}_0,
\,\tilde{\boldsymbol\Upsilon}_1,\,\tilde{\boldsymbol\Upsilon}_2,\,
\tilde{\boldsymbol\Upsilon}_3)$ to be positively polarized right orthonormal 
frames. Despite to this higher level of arbitrariness, here we can write the
relationships \mythetag{1.7} and the following relationships for
spinor frames:
$$
\xalignat 2
&\hskip -2em
\tilde{\boldsymbol\Psi}_i=\sum^4_{j=1}\hat{\goth S}^j_i
\ \boldsymbol\Psi_j,
&&\boldsymbol\Psi_i=\sum^4_{j=1}\hat{\goth T}^j_i
\ \tilde{\boldsymbol\Psi}_j.
\mytag{4.14}
\endxalignat
$$
Now $S^j_i$  and $\hat{\goth S}^j_i$ in \mythetag{1.7} and \mythetag{4.14}
are the components of arbitrary two $4\times 4$ matrices, while $T^j_i$ 
and $\hat{\goth T}^j_i$ are the components of their inverse matrices. 
Let's denote by $(U,\,\boldsymbol\vartheta^{\,1},\,\boldsymbol
\vartheta^{\,2},\,\boldsymbol\vartheta^{\,3},\,\boldsymbol\vartheta^{\,4})$
the dual frame for $(U,\,\boldsymbol\Psi_1,\,\boldsymbol\Psi_2,\,\boldsymbol
\Psi_3,\,\boldsymbol\Psi_4)$ and denote by $(\tilde U,\,\tilde{\boldsymbol
\vartheta}^{\,1},\,\tilde{\boldsymbol\vartheta}^{\,2},\,\tilde{\boldsymbol
\vartheta}^{\,3},\,\tilde{\boldsymbol\vartheta}^{\,4})$ the dual frame for 
$(\tilde U,\,\tilde{\boldsymbol\Psi}_1,\,\tilde{\boldsymbol\Psi}_2,\,
\tilde{\boldsymbol\Psi}_3,\,\tilde{\boldsymbol\Psi}_4)$. Then 
$$
\xalignat 2
&\hskip -2em
\tilde{\boldsymbol\vartheta}^{\,i}=\sum^4_{j=1}\hat{\goth T}^i_j
\ \boldsymbol\vartheta^{\,j},
&&\boldsymbol\vartheta^{\,i}=\sum^4_{j=1}\hat{\goth S}^i_j
\ \tilde{\boldsymbol\vartheta}^{\,j}.
\mytag{4.15}
\endxalignat
$$
Applying $\tau$ to $\boldsymbol\Psi_1,\,\boldsymbol\Psi_2,\,\boldsymbol
\Psi_3,\,\boldsymbol\Psi_4,\,\boldsymbol\vartheta^{\,1},\,\boldsymbol
\vartheta^{\,2},\,\boldsymbol\vartheta^{\,3},\,\boldsymbol\vartheta^{\,4},
\,\tilde{\boldsymbol\Psi}_1,\,\tilde{\boldsymbol\Psi}_2,\,
\tilde{\boldsymbol\Psi}_3,\,\tilde{\boldsymbol\Psi}_4,\,\tilde{\boldsymbol
\vartheta}^{\,1},\,\tilde{\boldsymbol\vartheta}^{\,2},\,\tilde{\boldsymbol
\vartheta}^{\,3}$,\linebreak and $\tilde{\boldsymbol\vartheta}^{\,4}$, we
get four frames in $D^{*\sssize\dagger}_p(M)$ and $D^{\sssize\dagger}_p(M)$:
$$
\xalignat 2
&\hskip -2em
\bPsi_i=\tau(\boldsymbol\Psi_i),
&&\overline{\boldsymbol\vartheta}\vphantom{\boldsymbol\vartheta}^{\,i}
=\tau(\boldsymbol\vartheta^{\,i}),
\mytag{4.16}\\
&\hskip -2em
\tilde{\bPsi}_i=\tau(\tilde{\boldsymbol\Psi}_i),
&&
\tilde{\overline{\boldsymbol\vartheta}}\vphantom{\vartheta}^{\,i}
=\tau(\tilde{\boldsymbol\vartheta}^{\,i}).
\mytag{4.17}
\endxalignat
$$
The frames $(U,\,\bPsi_1,\,\bPsi_2,\,\bPsi_3,\,\bPsi_4)$, $(\tilde U,\,
\tilde{\bPsi}\vphantom{\bPsi}_1,\,\tilde{\bPsi}\vphantom{\bPsi}_2,\,
\tilde{\bPsi}\vphantom{\bPsi}_3,\,\tilde{\bPsi}\vphantom{\bPsi}_4)$,
$(U,\,\overline{\boldsymbol\vartheta}\vphantom{\boldsymbol\vartheta}^{\,1},
\,\overline{\boldsymbol\vartheta}\vphantom{\boldsymbol\vartheta}^{\,2},\,
\overline{\boldsymbol\vartheta}\vphantom{\boldsymbol\vartheta}^{\,3},\,
\overline{\boldsymbol\vartheta}\vphantom{\boldsymbol\vartheta}^{\,4})$,
and $(\tilde U,\,\tilde{\overline{\boldsymbol
\vartheta}}\vphantom{\vartheta}^{\,1},\,\tilde{\overline{\boldsymbol
\vartheta}}\vphantom{\vartheta}^{\,2},\,\tilde{\overline{\boldsymbol
\vartheta}}\vphantom{\vartheta}^{\,3},\,\tilde{\overline{\boldsymbol
\vartheta}}\vphantom{\vartheta}^{\,4})$ are related to each other as
follows: 
$$
\xalignat 2
&\hskip -2em
\tilde{\bPsi}_i=\sum^4_{j=1}\overline{\hat{\goth S}^j_i}\ \bPsi_j,
&&\bPsi_i=\sum^4_{j=1}\overline{\hat{\goth T}^j_i}
\ \tilde{\bPsi}_j,
\mytag{4.18}\\
&\hskip -2em
\tilde{\overline{\boldsymbol\vartheta}}\vphantom{\vartheta}^{\,i}
=\sum^4_{j=1}\overline{\hat{\goth T}^i_j}\ \tilde{\boldsymbol
\vartheta}^{\,j},
&&\tilde{\boldsymbol\vartheta}^{\,i}=\sum^4_{j=1}
\overline{\hat{\goth S}^i_j}\ \tilde{\overline{\boldsymbol
\vartheta}}\vphantom{\vartheta}^{\,j}.
\mytag{4.19}
\endxalignat
$$
The formulas \mythetag{4.18} and \mythetag{4.19} are derived from
\mythetag{4.14} and \mythetag{4.15} by applying the relationships
\mythetag{4.16} and \mythetag{4.17}. And finally, we have the 
relationships
$$
\xalignat 2
&\hskip -2em
\tilde{\boldsymbol\eta}^i=\sum^3_{j=0}T^i_j\ \boldsymbol\eta^j,
&&\boldsymbol\eta^i=\sum^3_{j=0}S^i_j\ \tilde{\boldsymbol\eta}^j,
\mytag{4.20}
\endxalignat
$$
where $(U,\,\boldsymbol\eta^0,\,\boldsymbol\eta^1,\,\boldsymbol\eta^2,\,
\boldsymbol\eta^3)$ and $(\tilde U,\,\tilde{\boldsymbol\eta}^0,\,
\tilde{\boldsymbol\eta}^1,\,\tilde{\boldsymbol\eta}^2,\,\tilde{\boldsymbol
\eta}^3)$ are the frames dual to the frames $(U,\,\boldsymbol\Upsilon_0,\,
\boldsymbol\Upsilon_1,\,\boldsymbol\Upsilon_2,\,\boldsymbol\Upsilon_3)$
and $(\tilde U,\,\tilde{\boldsymbol\Upsilon}_0,\,\tilde{\boldsymbol
\Upsilon}_1,\,\tilde{\boldsymbol\Upsilon}_2,\,\tilde{\boldsymbol
\Upsilon}_3)$ respectively.\par
     Let's use the above frame vectors, covectors, spinors, and 
cospinors (including complex conjugate ones) in order to introduce 
the following tensor products:
$$
\align
&\hskip -2em
\aligned
\boldsymbol\Upsilon^{k_1\ldots\,k_n}_{h_1\ldots\,h_m}
&=\boldsymbol\Upsilon_{h_1}\otimes\ldots\otimes\boldsymbol\Upsilon_{h_m}
\otimes\boldsymbol\eta^{k_1}\otimes\ldots\otimes\boldsymbol\eta^{k_n},\\
\tilde{\boldsymbol\Upsilon}^{k_1\ldots\,k_n}_{h_1\ldots\,h_m}
&=\tilde{\boldsymbol\Upsilon}_{h_1}\otimes\ldots\otimes
\tilde{\boldsymbol\Upsilon}_{h_m}\otimes\tilde{\boldsymbol\eta}^{k_1}
\otimes\ldots\otimes\tilde{\boldsymbol\eta}^{k_n},
\endaligned
\mytag{4.21}\\
\vspace{2ex}
&\hskip -2em
\aligned
\boldsymbol\Psi^{j_1\ldots\,j_\beta}_{i_1\ldots\,i_\alpha}
&=\boldsymbol\Psi_{i_1}\otimes\ldots\otimes\boldsymbol\Psi_{i_\alpha}
\otimes\boldsymbol\vartheta^{\,j_1}\otimes\ldots\otimes
\boldsymbol\vartheta^{\,j_\beta},\\
\tilde{\boldsymbol\Psi}^{j_1\ldots\,j_\beta}_{i_1\ldots\,i_\alpha}
&=\tilde{\boldsymbol\Psi}_{i_1}\otimes\ldots\otimes
\tilde{\boldsymbol\Psi}_{i_\alpha}\otimes
\tilde{\boldsymbol\vartheta}^{\,j_1}\otimes\ldots\otimes
\tilde{\boldsymbol\vartheta}^{\,j_\beta},
\endaligned
\mytag{4.22}\\
\vspace{2ex}
&\hskip -2em
\aligned
\bPsi^{\bar j_1\ldots\,\bar j_\gamma}_{\bar i_1\ldots\,\bar i_\nu}
&=\bPsi_{\bar i_1}\otimes\ldots\otimes\bPsi_{\bar i_\nu}\otimes
\overline{\boldsymbol\vartheta}^{\,j_1}\otimes\ldots\otimes
\overline{\boldsymbol\vartheta}^{\,j_\gamma},\\
\tilde{\bPsi}\vphantom{\boldsymbol\Psi}^{\bar j_1\ldots\,
\bar j_\gamma}_{\bar i_1\ldots\,\bar i_\nu}
&=\tilde{\bPsi}_{\bar i_1}\otimes\ldots\otimes\tilde{\bPsi}_{\bar i_\nu}
\otimes\tilde{\overline{\boldsymbol\vartheta}}\vphantom{\vartheta}^{\,j_1}
\otimes\ldots\otimes
\tilde{\overline{\boldsymbol\vartheta}}\vphantom{\vartheta}^{\,j_\gamma}.
\endaligned
\mytag{4.23}
\endalign
$$
Then, using \mythetag{4.21}, \mythetag{4.22}, and \mythetag{4.23}, we
define other two tensor products:
$$
\gather
\hskip -2em
\boldsymbol\Psi^{j_1\ldots\,j_\beta\,\bar j_1\ldots\,\bar j_\gamma
\,k_1\ldots\, k_n}_{i_1\ldots\,i_\alpha\,\bar i_1\ldots\,\bar i_\nu
\,h_1\ldots\,h_m}
=\boldsymbol\Psi^{j_1\ldots\,j_\beta}_{i_1\ldots\,i_\alpha}
\otimes\bPsi^{\bar j_1\ldots\,\bar j_\gamma}_{\bar i_1\ldots\,
\bar i_\nu}\otimes\boldsymbol\Upsilon^{k_1\ldots\,k_n}_{h_1\ldots\,h_m},
\mytag{4.24}\\
\displaybreak
\hskip -2em
\tilde{\boldsymbol\Psi}^{j_1\ldots\,j_\beta\,\bar j_1\ldots\,\bar j_\gamma
\,k_1\ldots\, k_n}_{i_1\ldots\,i_\alpha\,\bar i_1\ldots\,\bar i_\nu\,
h_1\ldots\,h_m}=\tilde{\boldsymbol\Psi}^{j_1\ldots\,j_\beta}_{i_1\ldots
\,i_\alpha}\otimes\tilde{\bPsi}\vphantom{\boldsymbol\Psi}^{\bar j_1\ldots
\,\bar j_\gamma}_{\bar i_1\ldots\,\bar i_\nu}\otimes\tilde{\boldsymbol
\Upsilon}^{k_1\ldots\,k_n}_{h_1\ldots\,h_m}.
\mytag{4.25}
\endgather
$$
Both tensor products \mythetag{4.24} and \mythetag{4.25} are spin-tensorial
fields of the same type $(\alpha,\beta|\nu,\gamma|m,n)$. They are used in 
order to expand other spin-tensorial field of this type. If $\bold X$ is a
spin-tensorial field of the type $(\alpha,\beta|\nu,\gamma|m,n)$, then
$$
\bold X=\msum{4}\Sb i_1,\,\ldots,\,i_\alpha\\ j_1,\,\ldots,\,j_\beta
\endSb\msum{4}\Sb\bar i_1,\,\ldots,\,\bar i_\nu\\ \bar j_1,\,\ldots,\,
\bar j_\gamma\endSb\msum{3}\Sb h_1,\,\ldots,\,h_m\\ k_1,\,\ldots,\,k_n
\endSb
X^{i_1\ldots\,i_\alpha\bar i_1\ldots\,\bar i_\nu h_1\ldots\,h_m}_{j_1
\ldots\,j_\beta\bar j_1\ldots\,\bar j_\gamma k_1\ldots\, k_n}\ 
\boldsymbol\Psi^{j_1\ldots\,j_\beta\bar j_1\ldots\,\bar j_\gamma
k_1\ldots\, k_n}_{i_1\ldots\,i_\alpha\bar i_1\ldots\,\bar i_\nu
h_1\ldots\,h_m}.\quad
\mytag{4.26}
$$
The coefficients $X^{i_1\ldots\,i_\alpha\bar i_1\ldots\,\bar i_\nu
h_1\ldots\,h_m}_{j_1\ldots\,j_\beta\bar j_1\ldots\,\bar j_\gamma k_1
\ldots\, k_n}$ in \mythetag{4.26} are called the {\it component\/} of 
the field $\bold X$ in the pair of frames $(U,\,\boldsymbol\Psi_1,\,
\boldsymbol\Psi_2,\,\boldsymbol\Psi_3,\,\boldsymbol\Psi_4)$ and $(U,\,
\boldsymbol\Upsilon_0,\,\boldsymbol\Upsilon_1,\,\boldsymbol
\Upsilon_2,\,\boldsymbol\Upsilon_3)$. Similarly, the coefficients 
$\tilde X^{i_1\ldots\,i_\alpha\bar i_1\ldots\,\bar i_\nu h_1\ldots\,
h_m}_{j_1\ldots\,j_\beta\bar j_1\ldots\,\bar j_\gamma k_1\ldots\,
k_n}$ in the expansion 
$$
\bold X=\msum{4}\Sb i_1,\,\ldots,\,i_\alpha\\ j_1,\,\ldots,\,j_\beta
\endSb\msum{4}\Sb\bar i_1,\,\ldots,\,\bar i_\nu\\ \bar j_1,\,\ldots,\,
\bar j_\gamma\endSb\msum{3}\Sb h_1,\,\ldots,\,h_m\\ k_1,\,\ldots,\,k_n
\endSb
\tilde X^{i_1\ldots\,i_\alpha\bar i_1\ldots\,\bar i_\nu h_1\ldots\,
h_m}_{j_1\ldots\,j_\beta\bar j_1\ldots\,\bar j_\gamma k_1\ldots\, k_n}\ 
\tilde{\boldsymbol\Psi}^{j_1\ldots\,j_\beta\bar j_1\ldots\,\bar j_\gamma
k_1\ldots\, k_n}_{i_1\ldots\,i_\alpha\bar i_1\ldots\,\bar i_\nu h_1\ldots
\,h_m}\quad
\mytag{4.27}
$$
are called the components of the field $\bold X$ in the pair of frames 
$(\tilde U,\,\tilde{\boldsymbol\Psi}_1,\,\tilde{\boldsymbol\Psi}_2,\,
\tilde{\boldsymbol\Psi}_3,\,\tilde{\boldsymbol\Psi}_4)$ and $(\tilde U,
\,\tilde{\boldsymbol\eta}^0,\,\tilde{\boldsymbol\eta}^1,\,\tilde{\boldsymbol
\eta}^2,\,\tilde{\boldsymbol\eta}^3)$. Applying \mythetag{1.7}, 
\mythetag{4.14}, \mythetag{4.15}, \mythetag{4.18}, \mythetag{4.19}, and
\mythetag{4.20} to \mythetag{4.26} and \mythetag{4.27}, we derive the
following relationships
$$
\align
&\hskip -4em
\aligned
&\tilde X^{i_1\ldots\,i_\alpha\bar i_1\ldots\,\bar i_\nu
h_1\ldots\,h_m}_{j_1\ldots\,j_\beta\bar j_1\ldots\,\bar j_\gamma
k_1\ldots\,k_n}
=\dsize\msum{4}\Sb a_1,\,\ldots,\,a_\alpha\\ b_1,\,\ldots,\,b_\beta\endSb
\dsize\msum{4}\Sb \bar a_1,\,\ldots,\,\bar a_\nu\\ 
\bar b_1,\,\ldots,\,\bar b_\gamma\endSb
\dsize\msum{3}
\Sb c_1,\,\ldots,\,c_m\\ d_1,\,\ldots,\,d_n\endSb
\hat{\goth T}^{\,i_1}_{a_1}\ldots\,\hat{\goth T}^{\,i_\alpha}_{a_\alpha}
\,\times\\
&\kern 40pt
\times\,\hat{\goth S}^{b_1}_{j_1}\ldots\,\hat{\goth S}^{b_\beta}_{j_\beta}
\ \overline{\hat{\goth T}^{\,\bar i_1}_{\bar a_1}}
\ldots\,\overline{\hat{\goth T}^{\,\bar i_\nu}_{\bar a_\nu}}\ 
\ \overline{\hat{\goth S}^{\,\bar b_1}_{\bar j_1}}
\ldots\,\overline{\hat{\goth S}^{\,\bar b_\gamma}_{\bar j_\gamma}}\ 
T^{h_1}_{c_1}\ldots\,T^{h_m}_{c_m}\,\times\\
\vspace{1.5ex}
&\kern 60pt
\times\,S^{\,d_1}_{k_1}\ldots\,S^{\,d_n}_{k_n}\ 
X^{\,a_1\ldots\,a_\alpha\,\bar a_1\ldots\,\bar a_\nu\,
c_1\ldots\,c_m}_{\,b_1\ldots\,b_\beta\,\bar b_1\ldots\,\bar b_\gamma
\,d_1\ldots\,d_n},
\endaligned
\mytag{4.28}\\
\vspace{2ex}
&\hskip -4em
\aligned
&X^{i_1\ldots\,i_\alpha\bar i_1\ldots\,\bar i_\nu
h_1\ldots\,h_m}_{j_1\ldots\,j_\beta\bar j_1\ldots\,\bar j_\gamma
k_1\ldots\,k_n}
=\dsize\msum{4}\Sb a_1,\,\ldots,\,a_\alpha\\ b_1,\,\ldots,\,b_\beta\endSb
\dsize\msum{4}\Sb \bar a_1,\,\ldots,\,\bar a_\nu\\ 
\bar b_1,\,\ldots,\,\bar b_\gamma\endSb
\dsize\msum{3}
\Sb c_1,\,\ldots,\,c_m\\ d_1,\,\ldots,\,d_n\endSb
\hat{\goth S}^{\,i_1}_{a_1}\ldots\,\hat{\goth S}^{\,i_\alpha}_{a_\alpha}
\,\times\\
&\kern 40pt
\times\,\hat{\goth T}^{b_1}_{j_1}\ldots\,\hat{\goth T}^{b_\beta}_{j_\beta}
\ \overline{\hat{\goth S}^{\,\bar i_1}_{\bar a_1}}
\ldots\,\overline{\hat{\goth S}^{\,\bar i_\nu}_{\bar a_\nu}}\ 
\ \overline{\hat{\goth T}^{\,\bar b_1}_{\bar j_1}}
\ldots\,\overline{\hat{\goth T}^{\,\bar b_\gamma}_{\bar j_\gamma}}\ 
S^{h_1}_{c_1}\ldots\,S^{h_m}_{c_m}\,\times\\
\vspace{1.5ex}
&\kern 60pt
\times\,T^{\,d_1}_{k_1}\ldots\,T^{\,d_n}_{k_n}\ 
\tilde X^{\,a_1\ldots\,a_\alpha\,\bar a_1\ldots\,\bar a_\nu\,
c_1\ldots\,c_m}_{\,b_1\ldots\,b_\beta\,\bar b_1\ldots\,\bar b_\gamma
\,d_1\ldots\,d_n}.
\endaligned
\mytag{4.29}
\endalign
$$
The formulas \mythetag{4.28} and  \mythetag{4.29} represent the general
transformation rule for the components of spin-tensors in the case of
Dirac bundle $DM$. They are inverse to each other. Below we shall see
various special cases of them.\par
     Let's return back to the semilinear isomorphism of complex conjugation
$\tau$. From \mythetag{4.7} and \mythetag{4.8} we derive the following
relationships for $\tau$:
$$
\xalignat 2
&\hskip -2em
\tau(\boldsymbol\Upsilon_i)=\boldsymbol\Upsilon_i,
&&\tau(\boldsymbol\eta^i)=\boldsymbol\eta^i.
\mytag{4.30}
\endxalignat
$$ 
Now, if we combine \mythetag{4.30} with \mythetag{4.16} and \mythetag{4.17}
and if we remember the identity $\tau^2=\tau\compos\tau=\idop$, then we
derive the following formula:
$$
\tau(\bold X)=\msum{4}\Sb i_1,\,\ldots,\,i_\alpha\\ j_1,\,\ldots,
\,j_\beta\\ \bar i_1,\,\ldots,\,\bar i_\nu\\ \bar j_1,\,\ldots,\,
\bar j_\gamma\endSb\msum{3}\Sb h_1,\,\ldots,\,h_m\\ k_1,\,\ldots,
\,k_n\endSb
\overline{X^{\bar i_1\ldots\,\bar i_\alpha\,i_1\ldots\,i_\nu\,h_1\ldots\,
h_m}_{\bar j_1\ldots\,\bar j_\beta\,j_1\ldots\,j_\gamma\,k_1\ldots
\, k_n}}\ 
\boldsymbol\Psi^{j_1\ldots\,j_\gamma\,\bar j_1\ldots\,\bar j_\beta\,
k_1\ldots\, k_n}_{i_1\ldots\,i_\nu\,\bar i_1\ldots\,\bar i_\alpha\,
h_1\ldots\,h_m}.
\quad
\mytag{4.31}
$$
The formula \mythetag{4.31} means that the isomorphism $\tau$ acts upon
the components of the expansion \mythetag{4.26} as the complex conjugation 
exchanging barred and non-barred indices of them.\par
\head
5. Coordinate representation\\
of the basic spin-tensorial fields.
\endhead
    Let $(U,\,\boldsymbol\Psi_1,\,\boldsymbol\Psi_2)$ be an orthonormal 
frame of the chiral bundle $SM$. Then it induces three other orthonormal 
frames: $(U,\,\boldsymbol\vartheta^{\,1},\,\boldsymbol\vartheta^{\,2})$ 
in $S^*\!M$, $(U,\,\bPsi_1,\,\bPsi_2)$ in $S^{*\sssize\dagger}\!M$, and
$(U,\,\overline{\boldsymbol\vartheta}\vphantom{\boldsymbol
\vartheta}^{\,1},\overline{\boldsymbol\vartheta}\vphantom{\boldsymbol
\vartheta}^{\,2})$ in $S^{\sssize\dagger}\!M$. Due to the expansion 
\mythetag{3.1} two frames $(U,\,\boldsymbol\Psi_1,\,\boldsymbol\Psi_2)$
and $(U,\,\overline{\boldsymbol\vartheta}\vphantom{\boldsymbol
\vartheta}^{\,1},\overline{\boldsymbol\vartheta}\vphantom{\boldsymbol
\vartheta}^{\,2})$ compose a frame in $DM$. Let's denote
$$
\xalignat 2
&\hskip -2em
\boldsymbol\Psi_3=\overline{\boldsymbol\vartheta}
\vphantom{\boldsymbol\vartheta}^{\,1},
&&\boldsymbol\Psi_4=\overline{\boldsymbol\vartheta}
\vphantom{\boldsymbol\vartheta}^{\,2}.
\mytag{5.1}
\endxalignat
$$
\mydefinition{5.1}A frame $(U,\,\boldsymbol\Psi_1,\,\boldsymbol\Psi_2,
\,\boldsymbol\Psi_3,\,\boldsymbol\Psi_4)$ of the Dirac bundle $DM$ 
produced from some orthonormal frame $(U,\,\boldsymbol\Psi_1,\,
\boldsymbol\Psi_2)$ of the chiral bundle $SM$ by virtue of the formula 
\mythetag{5.1} is called a {\it canonically orthonormal chiral frame\/} 
of $DM$.
\enddefinition
     Let's consider the spin-metric tensor $\bold d$ introduced by the 
formula \mythetag{3.4}. In a canonically orthonormal chiral frame $(U,\,
\boldsymbol\Psi_1,\,\boldsymbol\Psi_2,\,\boldsymbol\Psi_3,\,\boldsymbol
\Psi_4)$ it is given by the matrix
$$
\hskip -2em
d_{ij}=
d(\boldsymbol\Psi_i,\boldsymbol\Psi_j)
=\Vmatrix 0 & 1 & 0 & 0\\-1 & 0 & 0 & 0\\
0 & 0 & 0 & -1\\0 & 0 & 1 & 0\endVmatrix.
\mytag{5.2}
$$
The matrix \mythetag{5.2} is a block-diagonal matrix composed of two
diagonal blocks. Its upper left diagonal block coincide with the matrix 
\mythetag{1.8} and its lower right diagonal block is given by the matrix
inverse to \mythetag{1.8}:
$$
\bd^{\,i\kern 0.5ptj}
=\bd(\overline{\boldsymbol\vartheta}\vphantom{\boldsymbol\vartheta}^{\,i},
\overline{\boldsymbol\vartheta}\vphantom{\boldsymbol\vartheta}^{\,j})
=\overline{d(\boldsymbol\vartheta^{\,i},\boldsymbol\vartheta^{\,j})}
=\Vmatrix 0 & -1\\1 & 0\endVmatrix.
$$ 
\mydefinition{5.2}A frame $(U,\,\boldsymbol\Psi_1,\,\boldsymbol\Psi_2,
\,\boldsymbol\Psi_3,\,\boldsymbol\Psi_4)$ of the Dirac bundle $DM$ 
is called an {\it orthonormal frame\/} if the spin-metric tensor $\bold d$
is represented by the matrix \mythetag{5.2} in this frame.
\enddefinition
     The chirality operator $\bold H$ is introduced by the formula
\mythetag{3.5}. It is easy to see that in a canonically orthonormal 
chiral frame it is represented by the matrix
$$
\hskip -2em
H^i_j=\Vmatrix 1 & 0 & 0 & 0\\0 & 1 & 0 & 0\\
0 & 0 & -1 & 0\\0 & 0 & 0 & -1\endVmatrix.
\mytag{5.3}
$$
\mydefinition{5.3} A frame $(U,\,\boldsymbol\Psi_1,\,\boldsymbol\Psi_2,
\,\boldsymbol\Psi_3,\,\boldsymbol\Psi_4)$ of the Dirac bundle $DM$ is
called a {\it chiral frame\/} if the chirality operator $\bold H$ is
given by the matrix \mythetag{5.3} in this frame.
\enddefinition
     The Hermitian spin-metric tensor $\bold D$ (it is also called 
the Dirac form) is introduced by the formulas \mythetag{3.7} and 
\mythetag{3.8}. In a canonically orthonormal chiral frame $(U,\,
\boldsymbol\Psi_1,\,\boldsymbol\Psi_2,\,\boldsymbol\Psi_3,\,\boldsymbol
\Psi_4)$ it is given by the matrix
$$
\hskip -2em
D_{i\bar j}=D(\boldsymbol\Psi_{\bar j},\boldsymbol\Psi_i)
=\Vmatrix 0 & 0 & 1 & 0\\0 & 0 & 0 & 1\\
1 & 0 & 0 & 0\\0 & 1 & 0 & 0\endVmatrix.
\mytag{5.4}
$$
\mydefinition{5.4} A frame $(U,\,\boldsymbol\Psi_1,\,\boldsymbol\Psi_2,
\,\boldsymbol\Psi_3,\,\boldsymbol\Psi_4)$ of the Dirac bundle $DM$ is 
called a {\it self-adjoint frame\/} if the Dirac form $\bold D$ is given 
by the matrix \mythetag{5.4} in this frame.
\enddefinition
    The following theorem links together the above four definitions.
\mytheorem{5.1} A frame $(U,\,\boldsymbol\Psi_1,\,\boldsymbol\Psi_2,
\,\boldsymbol\Psi_3,\,\boldsymbol\Psi_4)$ of the Dirac bundle $DM$
is a canonically orthonormal chiral frame if and only if it is 
orthonormal, chiral, and self-adjoint at the same time.
\endproclaim
     The theorem~\mythetheorem{5.1} shows that three basic spin-tensorial
fields $\bold d$, $\bold H$, and $\bold D$ describe completely the chiral 
expansion \mythetag{3.1} of the Dirac bundle $DM$.\par
\head
6. Geometrization of the extended group homomorphism.
\endhead
     Assume that we have two canonically orthonormal chiral frames 
of the Dirac bundle $(U,\,\boldsymbol\Psi_1,\,\boldsymbol\Psi_2,\,
\boldsymbol\Psi_3,\,\boldsymbol\Psi_4)$ and $(\tilde U,\,
\tilde{\boldsymbol\Psi}_1,\,\tilde{\boldsymbol\Psi}_2,\,\tilde{\boldsymbol
\Psi}_3,\,\tilde{\boldsymbol\Psi}_4)$ with overlapping domains. They are 
associated with the orthonormal frames $(U,\,\boldsymbol\Psi_1,\,
\boldsymbol\Psi_2)$ and $(\tilde U,\,
\tilde{\boldsymbol\Psi}_1,\,\tilde{\boldsymbol\Psi}_2)$ of the chiral
bundle $SM$, which in turn are associated with two positively polarized
right orthonormal frames $(U,\,\boldsymbol\Upsilon_0,\,\boldsymbol
\Upsilon_1,\,\boldsymbol\Upsilon_2,\,\boldsymbol\Upsilon_3)$ and 
$(\tilde U,\,\tilde{\boldsymbol\Upsilon}_0,\,\tilde{\boldsymbol
\Upsilon}_1,\,\tilde{\boldsymbol\Upsilon}_2,\,\tilde{\boldsymbol
\Upsilon}_3)$ of the tangent bundle $TM$. Thus we can write the complete
set of transition formulas \mythetag{1.7}, \mythetag{1.9}, \mythetag{4.14}, 
\mythetag{4.15}, \mythetag{4.18}, \mythetag{4.19}, and \mythetag{4.20}. 
However, the transition matrices $\hat{\goth S}$, $\hat{\goth T}$, $S$, 
and $T$ are not arbitrary $4\times 4$ matrices in this case. All these
matrices are determined by the only one $2\times 2$ matrix $\goth S\in
\MatGrSL(2,\Bbb C)$. This matrix $\goth S$ itself and its inverse matrix
$\goth T=\goth S^{-1}$ are explicitly present in the formulas
\mythetag{1.9}. The matrices $S$ and $T=S^{-1}$ in \mythetag{1.7} and
\mythetag{4.20} are produced from $\goth S$ and $\goth T$ by means
of the homomorphism \mythetag{1.1}:
$$
\xalignat 2
&\hskip -2em
S=\varphi(\goth S),
&&T=\varphi(\goth T).
\mytag{6.1}
\endxalignat
$$
The matrices $\hat{\goth S}$ and $\hat{\goth T}=\hat{\goth S}^{-1}$ are
produced from $\goth S$ and $\goth T$ in a more explicit way. They are
block-diagonal matrices constructed as follows:
$$
\xalignat 2
&\hskip -2em
\hat{\goth S}
=\Vmatrix\format\ \c\ &\,\c\ &\c\\
\goth S &\vrule height 8pt 
depth 4pt&0\\
\vspace{-8pt}
\kern -16pt\vbox{\hrule width 32pt}\kern -32pt&&\\
\vspace{-1pt}
0 &\vrule height 12pt depth 3pt
&\goth T^{\sssize\dagger}
\endVmatrix,
&&\hat{\goth T}
=\Vmatrix\format\ \c\ &\,\c\ &\c\\
\goth T &\vrule height 8pt 
depth 4pt&0\\
\vspace{-8pt}
\kern -15pt\vbox{\hrule width 33pt}\kern -33pt&&\\
\vspace{-1pt}
0 &\vrule height 12pt depth 3pt
&\goth S^{\sssize\dagger}
\endVmatrix.
\mytag{6.2}
\endxalignat
$$
Canonically orthonormal chiral frames of the Dirac bundle $DM$ are
naturally associated with positively polarized right orthonormal 
frames of the tangent bundle $TM$. Comparing \mythetag{6.2} with
\mythetag{2.19}, we see that transition matrices relating these two
types of frames form a group isomorphic to $\MatGrSL(2,\Bbb C)$
and the group $\MatGrSO^+(1,3,\Bbb R)$ respectively. The formulas
\mythetag{6.1} 
and \mythetag{6.2} mean that canonically orthonormal chiral frames
in $DM$ and positively polarized right orthonormal frames in $TM$
provide a geometrization of the upper line in the commutative
diagram \mythetag{2.40}.\par
    The coordinate representations of the basic spin-tensorial fields
$\bold d$, $\bold H$, and $\bold D$ are invariant when we change a
canonically orthonormal chiral frame for another such frame. Indeed, 
we have the relationships
$$
\gather
\hskip -2em
d_{ij}=\sum^4_{k=1}\sum^4_{q=1}\hat{\goth T}^k_i\,\hat{\goth T}^q_j\,
d_{kq},
\mytag{6.3}\\
\hskip -2em
H^i_j=\sum^4_{k=1}\sum^4_{q=1}\hat{\goth S}^i_k\,\hat{\goth T}^q_j\,
H^k_q,
\mytag{6.4}\\
\hskip -2em
D_{i\bar j}=\sum^4_{k=1}\sum^4_{q=1}\hat{\goth T}^k_i\ 
\overline{\hat{\goth T}^{\bar q}_{\bar j}}\ D_{k\bar q},
\mytag{6.5}
\endgather
$$
where $\hat{\goth S}$ and $\hat{\goth T}$ are block-diagonal matrices of
the form \mythetag{6.2}, while the components of $\bold d$, $\bold H$, 
and $\bold D$ are given by the matrices \mythetag{5.2}, \mythetag{5.3}, and
\mythetag{5.4} respectively. The formulas \mythetag{6.3}, \mythetag{6.4},
and \mythetag{6.5} can be verified by direct calculations.\par
     Comparing \mythetag{6.3}, \mythetag{6.4}, and \mythetag{6.5} with 
the general formulas \mythetag{4.28} and \mythetag{4.29}, we see that 
the spin-metric tensor $\bold d$ is a spin-tensorial field of the type
$(0,2|0,0|0,0)$, the chirality operator $\bold H$ is a spin-tensorial 
field of the type $(1,1|0,0|0,0)$, the Hermitian spin-metric $\bold D$
is a spin-tensorial field of the type $(0,1|0,1|0,0)$.\par
     In order to extend the above geometric interpretation of the upper
line of the commutative diagram \mythetag{2.40} to its lower line we need
to use the matrices \mythetag{2.37} as transition matrices and apply
them to some canonically orthonormal chiral frame $(U,\,\boldsymbol\Psi_1,
\,\boldsymbol\Psi_2,\,\boldsymbol\Psi_3,\,\boldsymbol\Psi_4)$. By setting
$\hat{\goth S}=\hat P$ in \mythetag{4.14} we get
$$
\xalignat 4
&\tilde{\boldsymbol\Psi}_1=\boldsymbol\Psi_3,
&&\tilde{\boldsymbol\Psi}_2=\boldsymbol\Psi_4,
&&\tilde{\boldsymbol\Psi}_3=\boldsymbol\Psi_1,
&&\tilde{\boldsymbol\Psi}_4=\boldsymbol\Psi_2.
\quad
\mytag{6.6}
\endxalignat
$$
We choose the plus sign in both formulas \mythetag{2.37} for the sake
of certainty. Since $\hat P^2=\bold 1$ (see \mythetag{2.36}), from
$\hat{\goth S}=\hat P$ we get $\hat{\goth T}=\hat{\goth S}^{-1}=
\hat P$. Then we derive
$$
\gather
\hskip -2em
\tilde d_{ij}=\sum^4_{k=1}\sum^4_{q=1}\hat P^k_i\,\hat P^q_j
\,d_{kq}=-d_{ij},
\mytag{6.7}\\
\hskip -2em
\tilde H^i_j=\sum^4_{k=1}\sum^4_{q=1}\hat P^i_k\,\hat P^q_j
\,H^k_q=-H^i_j,
\mytag{6.8}\\
\hskip -2em
\tilde D_{i\bar j}=\sum^4_{k=1}\sum^4_{q=1}\hat P^k_i\ 
\overline{\hat P^{\bar q}_{\bar j}}\ D_{k\bar q}
=D_{i\bar j}.
\mytag{6.9}
\endgather
$$
Like \mythetag{6.3}, \mythetag{6.4}, and \mythetag{6.5}, the formulas
\mythetag{6.7}, \mythetag{6.8}, and \mythetag{6.9} are easily derived 
by direct calculations.\par
     Note that the components of the chirality operator $\bold H$ change
their signs in \mythetag{6.8}. Therefore, the frame constructed by means
of the formulas \mythetag{6.6} is not a chiral frame, it is an {\it
anti-chiral frame}. \pagebreak The components of the spin-metric tensor 
$\bold d$
also change their signs. Hence, the frame \mythetag{6.6} is not an
orthonormal frame in the sense of the definition~\mythedefinition{5.2}.
It should be called an {\it anti-orthonormal frame}, though this is not
a commonly used term.
\mydefinition{6.1} A frame $(\tilde U,\,\tilde{\boldsymbol\Psi}_1,\,
\tilde{\boldsymbol\Psi}_2,\,\tilde{\boldsymbol\Psi}_3,\,\tilde{\boldsymbol
\Psi}_4)$ produced from some canonically orthonormal chiral frame $(U,\,
\boldsymbol\Psi_1,\,\boldsymbol\Psi_2,\,\boldsymbol\Psi_3,\,\boldsymbol
\Psi_4)$ by means of the formulas \mythetag{6.6} is called a {\it
$P$-reverse anti-chiral frame\/} of the Dirac bundle $DM$.
\enddefinition
\mytheorem{6.1} A frame $(\tilde U,\,\tilde{\boldsymbol\Psi}_1,\,
\tilde{\boldsymbol\Psi}_2,\,\tilde{\boldsymbol\Psi}_3,\,\tilde{\boldsymbol
\Psi}_4)$ is a $P$-reverse anti-chiral frame of the Dirac bundle $DM$ if
and only if the components of the basic spin-tensorial fields $\bold d$,
$\bold H$, and $\bold D$ in this frame are given by the formulas
$$
\xalignat 3
&\tilde d_{ij}=-d_{ij},
&&\tilde H^i_j=-H^i_j,
&&\tilde D_{i\bar j}=D_{i\bar j},
\endxalignat
$$
where $d_{ij}$, $H^i_j$, $D_{i\bar j}$ are taken from the matrices
\mythetag{5.2}, \mythetag{5.3}, and \mythetag{5.4} respectively.
\endproclaim
    Note that the spacial inversion matrix $\hat P$ from \mythetag{2.37}
is associated with the matrix $P$ in \mythetag{2.2} by means of the
formula \mythetag{2.32}. Therefore, each $P$-reverse anti-chiral frame 
of the Dirac bundle $DM$ is canonically associated with some positively
polarized left orthonormal frame in $TM$. This association yields a
partial geometrization of the group homomorphism \mythetag{2.39} forming
the lower line in the diagram \mythetag{2.40}. In order complete this
scheme of geometrization in the next step we consider the time inversion
matrix $\hat T$ from \mythetag{2.37}. By setting $\hat{\goth S}=\hat T$ 
in \mythetag{4.14} we get 
$$
\xalignat 4
&\tilde{\boldsymbol\Psi}_1=i\,\boldsymbol\Psi_3,
&&\tilde{\boldsymbol\Psi}_2=i\,\boldsymbol\Psi_4,
&&\tilde{\boldsymbol\Psi}_3=-i\,\boldsymbol\Psi_1,
&&\tilde{\boldsymbol\Psi}_4=-i\,\boldsymbol\Psi_2,
\qquad\quad
\mytag{6.10}
\endxalignat
$$
where $(U,\,\boldsymbol\Psi_1,\,\boldsymbol\Psi_2,\,\boldsymbol\Psi_3,
\,\boldsymbol\Psi_4)$ is some canonically orthonormal chiral frame of
$DM$. Taking $\hat{\goth S}=\hat T$, due to \mythetag{2.36} we get 
$\hat{\goth T}=\hat{\goth S}^{-1}=\hat T$. Then we derive
$$
\gather
\hskip -2em
\tilde d_{ij}=\sum^4_{k=1}\sum^4_{q=1}\hat T^k_i\,\hat T^q_j
\,d_{kq}=d_{ij},
\mytag{6.11}\\
\hskip -2em
\tilde H^i_j=\sum^4_{k=1}\sum^4_{q=1}\hat T^i_k\,\hat T^q_j
\,H^k_q=-H^i_j,
\mytag{6.12}\\
\hskip -2em
\tilde D_{i\bar j}=\sum^4_{k=1}\sum^4_{q=1}\hat T^k_i\ 
\overline{\hat T^{\bar q}_{\bar j}}\ D_{k\bar q}
=-D_{i\bar j}.
\mytag{6.13}
\endgather
$$
The formulas \mythetag{6.11}, \mythetag{6.12}, and \mythetag{6.13}
are analogous to \mythetag{6.7}, \mythetag{6.8}, and \mythetag{6.9}.
They are derived by direct calculations. In \mythetag{6.12} we see
that the components of the chirality operator $\bold H$ change their
signs. This means that the frame \mythetag{6.10}, like the frame 
\mythetag{6.6}, is an anti-chiral frame. However, unlike \mythetag{6.6},
it is an orthonormal frame in the sense of the 
definition~\mythedefinition{5.2} and it is not a self-adjoint frame 
in the sense of the definition~\mythedefinition{5.4}. Due to
\mythetag{6.13} it is an {\it anti-self-adjoint frame}.
\mydefinition{6.2} A frame $(\tilde U,\,\tilde{\boldsymbol\Psi}_1,\,
\tilde{\boldsymbol\Psi}_2,\,\tilde{\boldsymbol\Psi}_3,\,\tilde{\boldsymbol
\Psi}_4)$ produced from some canonically orthonormal chiral frame $(U,\,
\boldsymbol\Psi_1,\,\boldsymbol\Psi_2,\,\boldsymbol\Psi_3,\,\boldsymbol
\Psi_4)$ by means of the formulas \mythetag{6.10} is called a {\it
$T$-reverse anti-chiral frame\/} of the Dirac bundle $DM$.
\enddefinition
\mytheorem{6.2} A frame $(\tilde U,\,\tilde{\boldsymbol\Psi}_1,\,
\tilde{\boldsymbol\Psi}_2,\,\tilde{\boldsymbol\Psi}_3,\,\tilde{\boldsymbol
\Psi}_4)$ is a $T$-reverse anti-chiral frame of the Dirac bundle $DM$ if
and only if the components of the basic spin-tensorial fields $\bold d$,
$\bold H$, and $\bold D$ in this frame are given by the formulas
$$
\xalignat 3
&\tilde d_{ij}=d_{ij},
&&\tilde H^i_j=-H^i_j,
&&\tilde D_{i\bar j}=-D_{i\bar j},
\endxalignat
$$
where $d_{ij}$, $H^i_j$, $D_{i\bar j}$ are taken from the matrices
\mythetag{5.2}, \mythetag{5.3}, and \mythetag{5.4} respectively.
\endproclaim
    In the last step of our geometrization scheme we take the product 
$\hat Q=\hat P\cdot\hat T$. Due to our choice of positive signs in both
formulas \mythetag{2.37} we get
$$
\hskip -2em
\hat Q=\gamma_0\cdot\gamma_1\cdot\gamma_2\cdot\gamma_3
=\Vmatrix\format \l&\quad\r&\quad\r&\quad\r\\
i &0 &0 &0\\ 0 &i &0 &0\\ 0 &0 &-i &0\\ 0 &0 &0 &-i\endVmatrix.
\mytag{6.14}
$$
The product \mythetag{6.14} is taken from \mythetag{2.24}. By setting
$\hat{\goth S}=\hat Q$ in \mythetag{4.14} we get 
$$
\xalignat 4
&\tilde{\boldsymbol\Psi}_1=i\,\boldsymbol\Psi_1,
&&\tilde{\boldsymbol\Psi}_2=i\,\boldsymbol\Psi_2,
&&\tilde{\boldsymbol\Psi}_3=-i\,\boldsymbol\Psi_3,
&&\tilde{\boldsymbol\Psi}_4=-i\,\boldsymbol\Psi_4.
\qquad\quad
\mytag{6.15}
\endxalignat
$$
From \mythetag{2.38} we derive $\hat Q^2=-\bold 1$. Therefore, taking
$\hat{\goth S}=\hat Q$, we get $\hat{\goth T}=\hat{\goth S}^{-1}=-\hat Q$.
Then from \mythetag{4.28} we derive the relationships analogous to
\mythetag{6.11}, \mythetag{6.12}, and \mythetag{6.13}:
$$
\gather
\hskip -2em
\tilde d_{ij}=\sum^4_{k=1}\sum^4_{q=1}\hat Q^k_i\,\hat Q^q_j
\,d_{kq}=-d_{ij},
\mytag{6.16}\\
\hskip -2em
\tilde H^i_j=-\sum^4_{k=1}\sum^4_{q=1}\hat Q^i_k\,\hat Q^q_j
\,H^k_q=H^i_j,
\mytag{6.17}\\
\hskip -2em
\tilde D_{i\bar j}=\sum^4_{k=1}\sum^4_{q=1}\hat Q^k_i\ 
\overline{\hat Q^{\bar q}_{\bar j}}\ D_{k\bar q}
=-D_{i\bar j}.
\mytag{6.18}
\endgather
$$
The formulas \mythetag{6.16}, \mythetag{6.17}, and \mythetag{6.18} mean
that the frame \mythetag{6.15} is an anti-orthonormal, chiral, and
anti-self-adjoint frame of $DM$.
\mydefinition{6.3} A frame $(\tilde U,\,\tilde{\boldsymbol\Psi}_1,\,
\tilde{\boldsymbol\Psi}_2,\,\tilde{\boldsymbol\Psi}_3,\,\tilde{\boldsymbol
\Psi}_4)$ produced from some canonically orthonormal chiral frame $(U,\,
\boldsymbol\Psi_1,\,\boldsymbol\Psi_2,\,\boldsymbol\Psi_3,\,\boldsymbol
\Psi_4)$ by means of the formulas \mythetag{6.15} is called a {\it
$PT$-reverse chiral frame\/} of the Dirac bundle $DM$.
\enddefinition
\mytheorem{6.3} A frame $(\tilde U,\,\tilde{\boldsymbol\Psi}_1,\,
\tilde{\boldsymbol\Psi}_2,\,\tilde{\boldsymbol\Psi}_3,\,\tilde{\boldsymbol
\Psi}_4)$ is a $PT$-reverse chiral frame of the Dirac bundle $DM$ if
and only if the components of the basic spin-tensorial fields $\bold d$,
$\bold H$, and $\bold D$ in this frame are given by the formulas
$$
\xalignat 3
&\tilde d_{ij}=-d_{ij},
&&\tilde H^i_j=H^i_j,
&&\tilde D_{i\bar j}=-D_{i\bar j},
\endxalignat
$$
where $d_{ij}$, $H^i_j$, $D_{i\bar j}$ are taken from the matrices
\mythetag{5.2}, \mythetag{5.3}, and \mythetag{5.4} respectively.
\endproclaim
    Thus, the geometrization of the group homomorphism \mythetag{2.39}
is complete. The following diagram illustrates the frame association for
frames in $DM$ and $TM$:
$$
\hskip -2em
\aligned
&\boxit{Canonically orthonormal}{chiral frames}\to
\boxit{Positively polarized}{right orthonormal frames}\\
&\boxit{$P$-reverse}{anti-chiral frames}\to
\boxit{Positively polarized}{left orthonormal frames}\\
&\boxit{$T$-reverse}{anti-chiral frames}\to
\boxit{Negatively polarized}{right orthonormal frames}\\
&\boxit{$PT$-reverse}{chiral frames}\to
\boxit{Negatively polarized}{left orthonormal frames}
\endaligned
\mytag{6.19}
$$
Transition matrices relating frames in the right column of the
diagram \mythetag{6.19} form the group $\MatGrO(1,3,\Bbb R)$. For
frames in the left column their transition matrices form the
$4$-dimensional complex representation of the group $\GrPin(1,3,
\Bbb R)$.
\head
7. Spin-tensorial interpretation of the Dirac matrices.
\endhead
     Let's denote by $\gamma^{\,i}_{j\kern 0.2pt k}$ the components of 
the $k$-th Dirac matrix $\gamma_k$ and consider the matrix equality
\mythetag{2.22}. If we denote $\hat{\goth T}=\hat{\goth S}^{-1}$, then 
we write it as 
$$
\hskip -2em
\sum^3_{k=0}S^k_m\ \gamma^{\,i}_{j\kern 0.2pt k}
=\sum^4_{r=1}\sum^4_{s=1}
\hat\goth S^i_r\ \hat\goth T^s_j\ \gamma^{\,r}_{sm}.
\mytag{7.1}
$$
Using the inverse matrix $T=S^{-1}$, from \mythetag{7.1} we derive the
following equality:
$$
\hskip -2em
\gamma^{\,i}_{j\kern 0.2pt k}
=\sum^4_{r=1}\sum^4_{s=1}\sum^3_{m=1}
\hat\goth S^i_r\ \hat\goth T^s_j\ T^m_k\ \gamma^{\,r}_{sm}.
\mytag{7.2}
$$
The formula \mythetag{7.2} is a special case of the general 
transformation formula \mythetag{4.29}. It means that the components 
of all $\gamma$-matrices taken together define a spin-tensorial field 
of the type $(1,1|0,0|0,1)$. We denote it $\boldsymbol\gamma$. Here 
are the numeric values of $\gamma^{\,i}_{j\kern 0.2pt k}$ taken from
the matrices $\gamma_0$, $\gamma_1$, $\gamma_2$, and $\gamma_3$, in
\mythetag{2.25}, \mythetag{2.26}, \mythetag{2.27}, and \mythetag{2.28}:
$$
\xalignat 4
&\hskip -2em
\gamma^{\,1}_{1\,0}=0,
&&\gamma^{\,1}_{2\,0}=0,
&&\gamma^{\,1}_{3\,0}=1,
&&\gamma^{\,1}_{4\,0}=0,
\quad\\
&\hskip -2em
\gamma^{\,2}_{1\,0}=0,
&&\gamma^{\,2}_{2\,0}=0,
&&\gamma^{\,2}_{3\,0}=0,
&&\gamma^{\,2}_{4\,0}=1,
\quad\\
\vspace{-1.5ex}
&&&&&&&\mytag{7.3}\\
\vspace{-1.5ex}
&\hskip -2em
\gamma^{\,3}_{1\,0}=1,
&&\gamma^{\,3}_{2\,0}=0,
&&\gamma^{\,3}_{3\,0}=0,
&&\gamma^{\,3}_{4\,0}=0,
\quad\\
&\hskip -2em
\gamma^{\,4}_{1\,0}=0,
&&\gamma^{\,4}_{2\,0}=1,
&&\gamma^{\,4}_{3\,0}=0,
&&\gamma^{\,4}_{4\,0}=0,
\quad\\
\vspace{2ex}
&\hskip -2em
\gamma^{\,1}_{1\,1}=0,
&&\gamma^{\,1}_{2\,1}=0,
&&\gamma^{\,1}_{3\,1}=0,
&&\gamma^{\,1}_{4\,1}=1,
\quad\\
&\hskip -2em
\gamma^{\,2}_{1\,1}=0,
&&\gamma^{\,2}_{2\,1}=0,
&&\gamma^{\,2}_{3\,1}=1,
&&\gamma^{\,2}_{4\,1}=0,
\quad\\
\vspace{-1.5ex}
&&&&&&&\mytag{7.4}\\
\vspace{-1.5ex}
&\hskip -2em
\gamma^{\,3}_{1\,1}=0,
&&\gamma^{\,3}_{2\,1}=-1,
&&\gamma^{\,3}_{3\,1}=0,
&&\gamma^{\,3}_{4\,1}=0,
\quad\\
&\hskip -2em
\gamma^{\,4}_{1\,1}=-1,
&&\gamma^{\,4}_{2\,1}=0,
&&\gamma^{\,4}_{3\,1}=0,
&&\gamma^{\,4}_{4\,1}=0,
\quad\\
\displaybreak
&\hskip -2em
\gamma^{\,1}_{1\,2}=0,
&&\gamma^{\,1}_{2\,2}=0,
&&\gamma^{\,1}_{3\,2}=0,
&&\gamma^{\,1}_{4\,2}=-i,
\quad\\
&\hskip -2em
\gamma^{\,2}_{1\,2}=0,
&&\gamma^{\,2}_{2\,2}=0,
&&\gamma^{\,2}_{3\,2}=i,
&&\gamma^{\,2}_{4\,2}=0,
\quad\\
\vspace{-1.5ex}
&&&&&&&\mytag{7.5}\\
\vspace{-1.5ex}
&\hskip -2em
\gamma^{\,3}_{1\,2}=0,
&&\gamma^{\,3}_{2\,2}=i,
&&\gamma^{\,3}_{3\,2}=0,
&&\gamma^{\,3}_{4\,2}=0,
\quad\\
&\hskip -2em
\gamma^{\,4}_{1\,2}=-i,
&&\gamma^{\,4}_{2\,2}=0,
&&\gamma^{\,4}_{3\,2}=0,
&&\gamma^{\,4}_{4\,2}=0,
\quad\\
\vspace{2ex}
&\hskip -2em
\gamma^{\,1}_{1\,3}=0,
&&\gamma^{\,1}_{2\,3}=0,
&&\gamma^{\,1}_{3\,3}=1,
&&\gamma^{\,1}_{4\,3}=0,
\quad\\
&\hskip -2em
\gamma^{\,2}_{1\,3}=0,
&&\gamma^{\,2}_{2\,3}=0,
&&\gamma^{\,2}_{3\,3}=0,
&&\gamma^{\,2}_{4\,3}=-1,
\quad\\
\vspace{-1.5ex}
&&&&&&&\mytag{7.6}\\
\vspace{-1.5ex}
&\hskip -2em
\gamma^{\,3}_{1\,3}=-1,
&&\gamma^{\,3}_{2\,3}=0,
&&\gamma^{\,3}_{3\,3}=0,
&&\gamma^{\,3}_{4\,3}=0,
\quad\\
&\hskip -2em
\gamma^{\,4}_{1\,3}=0,
&&\gamma^{\,4}_{2\,3}=1,
&&\gamma^{\,4}_{3\,3}=0,
&&\gamma^{\,4}_{4\,3}=0.
\quad
\endxalignat
$$
In 
contrast to \mythetag{4.29}, we have no tilde in \mythetag{7.2}. 
Moreover, $\hat\goth S$ and $\hat\goth T$ are two mutually inverse 
matrices of the special form \mythetag{2.19}, while $T$ is produced 
from $\goth T$ by means of the group homomorphism \mythetag{1.1}. 
These features mean that $\gamma$-symbols given by \mythetag{7.3},
\mythetag{7.4}, \mythetag{7.5}, and \mythetag{7.6} should be ascribed
to canonically orthonormal chiral frames of $DM$ and to their associated
positively polarized right orthonormal frames in $TM$ (see the first line
in the diagram \mythetag{6.19}).\par
     Now let's proceed to the formulas \mythetag{2.32} and \mythetag{2.34}.
These two equalities can be easily transformed to the form similar to
\mythetag{7.2}:
$$
\align
&\hskip -2em
\gamma^{\,i}_{j\kern 0.2pt k}
=\sum^4_{r=1}\sum^4_{s=1}\sum^3_{m=1}
\hat P^i_r\ \hat P^s_j\ P^m_k\ \gamma^{\,r}_{sm},
\mytag{7.7}\\
&\hskip -2em
\gamma^{\,i}_{j\kern 0.2pt k}
=\sum^4_{r=1}\sum^4_{s=1}\sum^3_{m=1}
\hat T^i_r\ \hat T^s_j\ T^m_k\ \gamma^{\,r}_{sm}.
\mytag{7.8}
\endalign
$$
The formula \mythetag{7.7} is an analog of the formulas \mythetag{6.7},
\mythetag{6.8}, and \mythetag{6.9}, while \mythetag{7.8} is an analog 
of the formulas \mythetag{6.11}, \mythetag{6.12}, and \mythetag{6.13}.
From \mythetag{7.7} and \mythetag{7.8} one easily derives the following
formula for $\gamma$-symbols:
$$
\hskip -2em
\gamma^{\,i}_{j\kern 0.2pt k}
=-\sum^4_{r=1}\sum^4_{s=1}\sum^3_{m=1}
\hat Q^i_r\ \hat Q^s_j\ Q^m_k\ \gamma^{\,r}_{sm}.
\mytag{7.9}
$$
This formula \mythetag{7.9} is an analog of the formulas \mythetag{6.16},
\mythetag{6.17}, and \mythetag{6.18}. The matrix $\hat Q$ in it is taken
from the formula \mythetag{6.14}, while the matrix $Q$ is produced as the
product of the reflection matrices \mythetag{2.2}:
$$
Q=P\cdot T=T\cdot P=\Vmatrix
\format \l&\quad\r&\quad\r&\quad\r\\
-1 &0 &0 &0\\ 0 &-1 &0 &0\\ 0 &0 &-1 &0\\ 0 &0 &0 &-1\endVmatrix=-\bold 1.
$$
Due to the formulas \mythetag{7.7}, \mythetag{7.8}, and \mythetag{7.9}
the scope of the formulas \mythetag{7.3}, \mythetag{7.4}, \mythetag{7.5}, 
and \mythetag{7.6} can be extended so that we have the following theorem.
\mytheorem{7.1} Dirac's $\gamma$-symbols are the components of a 
spin-tensorial field of the type $(1,1|0,0|0,1)$ given by the formulas
\mythetag{7.3}, \mythetag{7.4}, \mythetag{7.5}, and \mythetag{7.6} in 
any frame pair specified in the diagram \mythetag{6.19}.
\endproclaim
\head
8. some conclusions.
\endhead
    The $P$ and $T$ operators are introduced in special relativity in
order to describe the spacial and time inversion operations for wave
functions of elementary particles:
$$
\align
&\hskip -2em
P\!:\quad\psi(t,x,y,z)\ \to\ \psi(t,-x,-y,-z),
\mytag{8.1}\\
&\hskip -2em
T\!:\quad\psi(t,x,y,z)\ \to\ \psi(-t,x,y,z)
\mytag{8.2}
\endalign
$$
(see \mycite{7} and \mycite{8} for details). However, in general 
relativity the coordinate transformations \mythetag{8.1} and
\mythetag{8.2} are not permitted, provided the space-time manifold
$M$ and its metric $\bold g$ are fixed. For this reason here $P$ 
and $T$ transformations are interpreted not as actual operators, 
but as frame transformations only. It seems to me, that in order 
to treat $P$ and $T$ as actual operators (as actual symmetries of
the Nature) one should add some transformations of $M$ and $\bold g$
performed simultaneously with the transformations \mythetag{8.1} and
\mythetag{8.2}.

\enddocument
\end